\documentclass[12pt,titlepage]{amsart}

\usepackage[T1]{fontenc} 
\usepackage[utf8]{inputenc} 
\usepackage{longtable}
\usepackage{amssymb}
\usepackage{amsthm}
\usepackage{amsmath}
\usepackage{wasysym}
\usepackage{dcpic, pictex}
\usepackage{tikz-cd}
\usepackage{bbm}
\usepackage{enumerate}
\usepackage{amsfonts}
\usepackage{amsthm}
\usepackage{amsmath}
\usepackage{times}
\usepackage{amscd}
\usepackage[mathscr]{eucal}
\usepackage{indentfirst}
\usepackage{graphicx}
\usepackage{graphics}
\usepackage{pict2e}
\usepackage{epic}
\usepackage{epstopdf} 
\usepackage{verbatim}
\usepackage{cancel}
\usepackage[foot]{amsaddr}
\usepackage{lipsum}

\usepackage{
amsfonts, 
latexsym, 
enumerate,
}
\usepackage[english]{babel}
\allowdisplaybreaks
\makeindex
\newtheorem{theorem}{Theorem}[section]
\newtheorem{lemma}[theorem]{Lemma}

\newtheorem{definition}[theorem]{Definition}

\newtheorem{remark}[theorem]{Remark}
\newtheorem{corollary}[theorem]{Corollary}
\usepackage[a4paper,top=2.9cm,bottom=2.9cm,left=2.1cm,right=2.1cm]{geometry}
\numberwithin{equation}{section}
\allowdisplaybreaks
\newcommand\numberthis{\addtocounter{equation}{1}\tag{\theequation}}
\usepackage{ bbold }
\usepackage{xcolor}

\begin{document}

\title{FLUCTUATIONS FOR LINEAR EIGENVALUE STATISTICS OF SAMPLE COVARIANCE MATRICES}

\author{Giorgio Cipolloni$^1$, L\text{\'{a}}szl\text{\'{o}} Erd\text{\H{o}}s$^1$}
\thanks{$^1$IST Austria, Am Campus 1, 3400 Klosterneuburg, Austria}
\thanks{E-mail addresses: gcipollo@ist.ac.at, lerdos@ist.ac.at}
\thanks{This project has received funding from the European Union's Horizon 2020 research and innovation programme \text{ } \text{ } \text{ } \text{ } \text{ $\,$}  under the Marie Sk\l odowska-Curie Grant Agreement No. 665385}
\thanks{Partially supported by the ERC Advanced grant RANMAT No. 338804.}

\begin{abstract}
We prove a central limit theorem for the difference of linear eigenvalue statistics of a sample covariance matrix $\widetilde{W}$ and its minor $W$. We find that the fluctuation of this difference is much smaller than those of the individual linear statistics, as a consequence of the strong correlation between the eigenvalues of $\widetilde{W}$ and $W$. Our result identifies the fluctuation of the spatial derivative of the approximate Gaussian field in the recent paper by Dumitru and Paquette. Unlike in a similar result for Wigner matrices, for sample covariance matrices the fluctuation may entirely vanish.
\end{abstract} 

{\let\newpage\relax\maketitle}

\section{Introduction}
We consider sample covariance matrices of the form $\widetilde{W}=\widetilde{X}^*\widetilde{X}$, where the entries of the $M\times N$ matrix $\widetilde{X}$ are i.i.d. random variables with mean zero and variance $\frac{1}{\sqrt{MN}}$. In the Gaussian case this ensemble was introduced by Wishart \cite{wishart1928generalised}. Besides Wigner matrices, this is the oldest and the most studied family of random matrices.

Let $\lambda_1,\dots,\lambda_N$ be the eigenvalues of $\widetilde{W}=\widetilde{X}^*\widetilde{X}$, then the empirical distribution $\frac{1}{N}\sum_{i=1}^N\delta_{\lambda_i}$ converges in probability to the Marchenko-Pastur distribution \cite{marvcenko1967distribution}. This asymptotics can be refined by examining the centered linear statistics \begin{equation}
\label{cent}
\text{Tr}f(\widetilde{W})-\mathbb{E}\text{Tr}\,f(\widetilde{W})=\sum_{i=1}^N\left[f(\lambda_i)-\mathbb{E}f(\lambda_i)\right],
\end{equation} with a sufficiently smooth function $f$, which has been shown to have Gaussian fluctuation (see e.g. \cite{bai2008clt}, \cite{jonsson1982some}, \cite{shcherbina2011central}). Notice that \eqref{cent} does not carry the usual $\frac{1}{\sqrt{N}}$ normalization of the conventional central limit theorem. In particular this result indicates a very strong correlation between eigenvalues. Apart from understanding an interesting mathematical phenomenon, the asymptotic properties of centered linear statistics for sample covariance matrices also have potential applications \cite{rao2008statistical}.

All the previously cited works on the centered linear statistics of a sample covariance matrix $\widetilde{W}$ concern the study of a single random matrix. The recent paper of Dumitru and Paquette \cite{dumitriu2018spectra} considers the joint eigenvalue fluctuations of a sample covariance matrix and its minors, by picking submatrices whose dimensions differ macroscopically. They show that their centered linear eigenvalue statistics converge to spatial averages of a two dimensional Gaussian free field. Similar results for Wigner matrices have been achieved earlier in \cite{borodin2014clt}.

In the current work we study this phenomenon for submatrices whose dimensions differ only by one. This requires a detailed analysis on the local spectral scale while \cite{dumitriu2018spectra} concerns only the global scale. In particular, we prove a central limit theorem (CLT) for the difference of linear eigenvalue statistics of a sample covariance matrix $\widetilde{W}=\widetilde{X}^*\widetilde{X}$ and its minor $W=X^*X$, obtained by deleting the first row and column. This difference fluctuates on a scale $N^{-\frac{1}{2}}$, which is much smaller than the order one fluctuations scale of the individual linear statistics, demonstrating a strong correlation between  the eigenvalues of $\widetilde{W}$ and its minor $W$. The statistical interpretation of our result is that changing the sample size by one in a statistical data has very little influence on the fluctuations of the linear eigenvalue statistics. Motivated by Gorin and Zhang \cite{gorinzhang}, another interpretation is that we prove a CLT for the \textit{spatial derivative} of the approximate Gaussian field in \cite{dumitriu2018spectra}.

This result extends a CLT, proved in \cite{erdHos2017fluctuations} for Wigner matrices, to sample covariance random matrices, with the difference that in this latter case it is also possible not to have random fluctuations at all, see Remark \ref{nofluc} in Section 2.

In the proof of the CLT for sample covariance matrices there are two main differences compared to the proof given in \cite{erdHos2017fluctuations} for the Wigner case.  Firstly, we have to handle the singularity of the Marchenko-Pastur law at zero, which also gives an additional contribution to the leading order term of \eqref{mainf}. Secondly, the entries of the matrix $\widetilde{W}=\widetilde{X}^*\widetilde{X}$ are not independent and the analogy occurs on the level of $\widetilde{X}$. Besides linearizing the problem and using recent local laws for Gram matrices \cite{alt2017local, bloemendal2014isotropic}, we need to approximate sums of the form $\sum_{ij} G_{ij}G_{ji}'$ and $\sum_{ij} G_{ij}G_{ij}'$ where $G$ and $G'$ are the resolvents of $XX^*$ at two different spectral parameter. While the first sum is tracial, the second one is not and thus cannot be directly analyzed by existing local laws: we need to derive a novel self-consistent equation for it.

\section*{Notation}
We introduce some notation we use throughout the paper. For positive quantities $f,g$, we write $f\lesssim g$ if $f\le Cg$, for some $C>0$ which depends only on the parameter $\phi$ defined in \eqref{limitN}. Similarly, we define $f\gtrsim g$. For any $\alpha,\beta>0$, with $\alpha \asymp \beta$ we denote that there exists two $\phi$ independent constants $r_*, r^*>0$ such that $r_* \beta \le \alpha \le r^*\beta$.

\section*{Acknowledgments}
The authors are grateful to Dominik Schr\text{\"{o}}der for valuable discussions. We also thank the referees for many useful comments and  for pointing out a computation mistake in the first version of the paper.

\section{Main Results}
All along the paper we will refer to the $N\times N$ matrix with $\widetilde{W}=\widetilde{X}^*\widetilde{X}$ and to the $(N-1)\times(N-1)$ matrix obtained after removing its first row and column with $W=X^*X$, where $X$ is the matrix obtained by $\widetilde{X}$ after removing its first column. It may look unconventional, but we chose to put the tilde on the original matrix $\widetilde{W}$ and no tilde on the minor $W$ in order to simplify formulas.

\begin{remark}
\label{ind}
We follow the convention that Latin letters $i\in\{1,\dots,M\}$ denote the rows of the matrix $\widetilde X$ and Greek letters $\mu\in\{1,\dots,N\}$ its columns.
\end{remark}

Let $\widetilde{X}$ be an $M\times N$ matrix whose entries $\widetilde{X}_{i\mu}$ are i.i.d. complex valued random variables satisfying: \begin{equation}
\label{var}
\mathbb{E}\widetilde{X}_{i\mu}=0, \,\, \mathbb{E}|\widetilde{X}_{i\mu}|^2=\frac{1}{\sqrt{MN}},\,\,\,\,\, 1\le i\le M,\,1\le \mu\le N.
\end{equation} Furthermore, for any $p\in\mathbb{N}$ there exists a constant $C_p>0$ such that \begin{equation}
\label{mom}
\mathbb{E}\left|(NM)^\frac{1}{4}\widetilde{X}_{i\mu}\right|^p\le C_p,\,\,\,\,\, 1\le i\le M,\,1\le \mu\le N.
\end{equation} We assume that $M$ and $N$ are comparable, i.e. there exist $N$-independent constants $c_1,c_2>0$ such that \begin{equation}
\label{limitN}
c_1\le \phi:=\frac{M}{N}\le c_2.
\end{equation}

For fixed $\phi$ and large $N$ the empirical distribution of the eigenvalues of the $N\times N$ matrix $\widetilde{W}=\widetilde{X}^*\widetilde{X}$ is given by the Marchenko-Pastur law \cite{marvcenko1967distribution}: \begin{equation}
\label{mplaw}
\rho_\phi(dx)=\rho_\phi(x)dx+(1-\phi)_+\delta(dx),\,\,\text{with}\,\,\,\rho_\phi(x):=\frac{\sqrt{\phi}}{2\pi} \sqrt{\frac{[(x-\gamma_-)(\gamma_+-x)]_+}{x^2}},
\end{equation} where we defined \[\gamma_\pm:=\sqrt{\phi}+\frac{1}{\sqrt{\phi}}\pm 2\] to be the edges of the limiting spectrum. The Stieltjes transform of $\rho_\phi(dx)$ is \begin{equation}
\label{st}
m_\phi(z):=\int_\mathbb{R} \frac{\rho_\phi(dx)}{x-z}=\frac{\phi^{1/2}-\phi^{-1/2}-z+i\sqrt{(z-\gamma_-)(\gamma_+-z)}}{2\phi^{-1/2}z},
\end{equation} where the square root is chosen so that $m_\phi$ is holomorphic in the complex upper half plane $\mathbb{H}$ and satisfies $m_\phi(z)\to 0$ as $z\to \infty$. The function $m_\phi=m_\phi(z)$ may also be characterized as the unique solution of the equation \begin{equation}
\label{st1}
m_\phi+\frac{1}{z+z\phi^{-1/2}m_\phi-(\phi^{1/2}-\phi^{-1/2})}=0
\end{equation} satisfying $\Im m_\phi(z)>0$ for $\Im z>0$. Our main result is the following:


\begin{theorem}
\label{maintheo}
Let $d_*>0$ and $\widetilde{W}=\widetilde{X}^*\widetilde{X}$, with $\widetilde{X}$ an $M\times N$ matrix whose i.i.d. entries satisfy \eqref{var} and \eqref{mom}. Furthermore we assume \eqref{limitN} and that either $\phi=1$ or $|\phi-1|\ge d_*$. Let $\sigma_2:=\sqrt{MN}\mathbb{E}\widetilde{X}_{i \mu}^2$ and $\sigma_4:=MN\mathbb{E}|\widetilde{X}_{i\mu}|^4$ and assume that $\sigma_2$ and $\sigma_4$ are $N$-independent . Moreover, let  $f\in H_0^2([\gamma_- -\delta_*,\gamma_+ +\delta_*])$, for a small $\delta_*=\delta_*(d_*)$ such that $\gamma_--\delta_*>0$ for $|\phi-1|\ge d_*$, be some real valued function in the $H_0^2$-Sobolev space. Then the random variable 
\begin{equation}
\label{mainf}
 f_N:=\text{Tr}f(\widetilde{W})-\text{Tr}f(W)
 \end{equation}
converges in probability to the constant 
\begin{equation}
\label{lott}
\Omega_f:=\int_{\gamma_-}^{\gamma_+}f(x) \frac{\sqrt{\phi}}{4\pi^2 x\rho_\phi(x)}\left(1+\frac{\sqrt{\phi}-\frac{1}{\sqrt{\phi}}}{x}\right)\,dx
\end{equation} 
for $|\phi-1|\ge d_*$, and to
\[
\Omega_f:=\int_0^4\frac{f(x)}{4\pi^2 x\rho_1(x)}\,dx+\frac{f(0)}{2}.
\] for $\phi=1$.
More precisely, for any fixed $\epsilon>0$, 
\[\mathbb{E}f_N=\Omega_f+\mathcal{O}\left(N^{-\frac{2}{3}+\epsilon}\right)\] 
and $f_N$ fluctuates on a scale $N^{-\frac{1}{2}}$, i.e. 
\[
\mathbb{E}\left(\sqrt{N}(f_N-\Omega_f)\right)^2=V_f+\mathcal{O}\left(N^{-\frac{1}{6}+\epsilon}\right).
\]
The limit variance $V_f$ can be computed explicitly: \begin{equation}
\label{var2}
V_f:=V_{f,1}+(\sigma_4-1)V_{f,2}+|\sigma_2|^2V_{\sigma_2},
\end{equation} with \[V_{f,1}=\int_{\gamma_-}^{\gamma_+} f'(x)^2x\rho_\phi(x)\phi^{-\frac{1}{2}}\,dx-\left(\int_{\gamma_-}^{\gamma_+} f'(x)x\rho_\phi(x)\phi^{-\frac{1}{2}}\, dx\right)^2,\]  \[V_{f,2}=\left(\int_{\gamma_-}^{\gamma_+} f'(x)x\rho_\phi(x)\phi^{-\frac{1}{2}}\, dx\right)^2,\] where $\rho_\phi(x)$ is the density of the Marchenko-Pastur law \eqref{mplaw}, and $V_{\sigma_2}$ defined as in \eqref{G} if $|\sigma_2|<1$ and $V_{\sigma_2}:=V_{f,1}$ if $|\sigma_2|=1$.

Furthermore, \[\sqrt{N}(f_N-\Omega_f)\Rightarrow \Delta_f,\] where $\Delta_f$ is a centered Gaussian random variable of variance $V_f$ and $"\Rightarrow"$ denotes the convergence in distribution. Finally, any fixed moment converges at least at a rate $\mathcal{O}\left(N^{-\frac{1}{6}+\epsilon}\right)$ to the corresponding Gaussian moments.
\end{theorem}

\begin{remark}
The non-negativity of $V_{f,1}$ follows by applying Schwarz inequality and by using that $\int_{\gamma_-}^{\gamma_+} x\rho_\phi(x)\phi^{-\frac{1}{2}}\,dx=1$.
\end{remark}

\begin{remark}
\label{nofluc}
One can easily check that the variance $V_f$ is zero if and only if $\sigma_2=0,\,\sigma_4=1$ and $f'(x)\equiv 1$. This is the case, for example, when the entries of $\widetilde{X}$ are i.i.d complex Bernoulli random variables, i.e. the distribution of each $\widetilde{X}_{i\mu}$ is $(MN)^{-\frac{1}{4}}e^{iU}$, with $U$ a uniform random variable in $[0,2\pi]$. In particular, since the entries of $\widetilde{X}$ have modulus $(MN)^{-\frac{1}{4}}$, the difference of the traces of $\widetilde{W}$ and $W$ is deterministic: \[\text{Tr}f(\widetilde{W})-\text{Tr}f(W)=\text{Tr} \widetilde{W}-\text{Tr}W= \bf{x}^*\bf{x}=\sqrt{\phi},\] where $\bf{x}$ is the first column of $\widetilde{X}$. The possibility of $V_f=0$ is a fundamental difference compared to the Wigner case in \cite{erdHos2017fluctuations} where the analogous quantity always had a non trivial fluctuation.
\end{remark}

\begin{remark}
We stated our result in Theorem \ref{maintheo} for the matrix $\widetilde{X}^*\widetilde{X}$, but it obviously holds for $\widetilde{X}\widetilde{X}^*$ as well. Indeed all computations and results remain valid after the swapping: $\widetilde{X}\leftrightarrow \widetilde{X}^*$, $M\leftrightarrow N$ and $\phi\leftrightarrow\phi^{-1}$. The empirical distribution of the eigenvalues of $\widetilde{X}\widetilde{X}^*$ is asymptotically  $\rho_{\phi^{-1}}(dx)$, whose Stieltjes transform is \begin{equation}
\label{change1}
m_{\phi^{-1}}(z)=\frac{1}{\phi}\left(m_\phi(z)+\frac{1-\phi}{z}\right).
\end{equation}
\end{remark}

\begin{remark}
Notice that in the statement of Theorem \ref{maintheo} we assumed that $\widetilde{X}$ is either a square matrix, $\phi=1$, or a proper rectangular matrix, $|\phi-1|>d_*$. The reason is that to prove Theorem \ref{maintheo} we use optimal local laws for all $z\in \mathbb{H}$ which are available in these cases only (see \cite{alt2017local}).
\end{remark}


\section{Preliminaries}
Our main result pertains to the matrix $\widetilde{X}^*\widetilde{X}$, but in the proof we will also need the matrix $\widetilde{X}\widetilde{X}^*$, so for each $z\in\mathbb{H}$ we define both resolvents \begin{equation}
\label{res}
\widetilde{R}(z):=(\widetilde{X}^*\widetilde{X}-z)^{-1},\,\,\, \widetilde{G}(z):=(\widetilde{X}\widetilde{X}^*-z)^{-1}.
\end{equation} 

Next, we define the $M\times (N-1)$ matrix $X$
as the matrix $\widetilde{X}$ after removing its first column, which we denote by ${\bf x}$, i.e. $\widetilde{X}=[{\bf x}|X]$. Moreover, for $\mu,\nu \notin \{1\}$, we define the resolvent entries \[
R_{\mu\nu}(z):=\left(X^*X-z\right)_{\mu\nu}^{-1},\,\,\, G_{ij}^{[T]}(z):=\left(XX^*-z\right)_{ij}^{-1}.
\]

\begin{remark}
In the following sections, without loss of generality, we will always assume that $\phi\ge 1$, i.e. $M\ge N$. Indeed, if $\phi\le 1-d_*$ then the proof proceeds exactly in the same way having in mind that $m_{\phi^{-1}}$ and $m_{\phi}$ are related by \eqref{change1}.
\end{remark}

Since $\phi\ge 1$ and the spectrum of $\widetilde{X}\widetilde{X}^*$ is equal to the spectrum of $\widetilde{X}^*\widetilde{X}$ plus $M-N$ zero eigenvalues, we have  \begin{equation}
\label{trace}
\phi\frac{1}{M}\text{Tr}\widetilde{G}=\frac{1}{ N}\text{Tr}\widetilde{R}+\frac{1-\phi}{z}
\end{equation} and that \begin{equation}
\label{tracem}
\text{Tr}R-\text{Tr}G=\frac{M-(N-1)}{z}.
\end{equation} Furthermore, setting $\eta=\Im z>0$, we have the \textit{Ward identity} \begin{equation}
\label{ward}
\sum_{j=1}^M\left|G_{ij}(z)\right|^2=\frac{1}{\eta}\Im G_{ii}(z).
\end{equation}

Finally, we record some properties of the Stieltjes transform defined in \eqref{st} in the following lemma, which will be proved in Appendix A.

\begin{lemma}
\label{sbound}
There exist positive constants $c,\tilde{c}, \hat{c}$ such that for any $\phi\ge 1$ and for each $z=x+i\eta\in\mathbb{H}$ such that $|z-\sqrt{\phi}|\le 10$ we have the following bounds

\begin{align}
c\le\left|\frac{z}{\sqrt{\phi}}m_\phi(z)^2\right|&\le 1-\tilde{c}\eta, \label{goodbound}\\
\left|m_\phi(z)'\right|&\le \frac{\hat{c}\sqrt{\phi}}{|z|\sqrt{\kappa_x+\eta}}, \label{boundderivative}\\
\left|1-z\phi^{-\frac{1}{2}}m_\phi(z)^2\right|&\asymp \frac{\phi^\frac{1}{4}}{|z|^\frac{1}{2}}\sqrt{\kappa_x+\eta}, \label{finaleq},
\end{align} where $\kappa_x:=\min\{|\gamma_+-x|,|\gamma_--x|\}$. For $c_1\le\phi\le 1-d_*$ the same bounds hold for $z=x+i\eta\in\mathbb{H}$ such that $\gamma_--\delta_*\le x \le 10$, with $\delta_*=\delta_*(d_*)$ such that  $\gamma_--\delta_*>0$.
\end{lemma}

In Lemma \ref{sbound} we explicitly wrote the $\phi$-dependence in the bounds since they hold uniformly in $\phi$. But all along the proof of Theorem \ref{maintheo} we will omit the explicit dependence on $\phi$, since we work under the assumption $c_1\le \phi \le c_2$ (see \eqref{limitN}).

\section{Mean and  variance  computation}
In this section we prove Theorem \ref{maintheo} in the sense of mean and variance. We recall that with ${\bf x}$ we denote the first column of $\widetilde{X}$. To study $f_N=\text{Tr} f(\widetilde{W})-\text{Tr}f(W)$, with $\widetilde{W}=\widetilde{X}^*\widetilde{X}$ and $W=X^*X$, we consider the quantity \begin{equation}
\label{delta}
\Delta_N(z):=\text{Tr}\widetilde{R}(z)-\text{Tr}R(z),\,\,\, z\in\mathbb{H}.
\end{equation}

Clearly $\widetilde{X}\widetilde{X}^*$ is a rank-one perturbation of the matrix $XX^*$, hence to compute $\widetilde{G}(z)$ we use the following lemma whose proof is a direct calculation.

\begin{lemma}
\label{expres}
Let $A$ be an $M\times M$ matrix with $\Im A<0$ and $h\in\mathbb{C}^M$ a column vector, then \[ \frac{1}{A+hh^*}=\frac{1}{A}-\frac{1}{1+\left\langle h,\frac{1}{A}h\right\rangle}\cdot \frac{1}{A} hh^*\frac{1}{A}.\]
\end{lemma}

We now find an explicit formula for $\Delta_N(z)$. Using \eqref{trace}, \eqref{tracem} and  \eqref{delta}
 we get
  \[\Delta_N(z)=\text{Tr}\widetilde{G}(z)-\text{Tr}R(z)-\frac{N(1-\phi)}{z}=\text{Tr}\frac{1}{XX^*+{\bf xx}^*-z}-\text{Tr}G(z)-\frac{1}{z}.\] Using Lemma \ref{expres} for the first term in the right-hand side, we conclude that \begin{equation}
\label{delta12}
\Delta_N(z)=-\frac{\left\langle {\bf x}, G^2(z){\bf x}\right\rangle}{1+\left\langle {\bf x},G(z){\bf x}\right\rangle}-\frac{1}{z}.
\end{equation}

We introduce a commonly used notion of high probability bound.

\begin{definition}
\label{sd}
If \[ X=\left(X^{(N)}(u)|N\in\mathbb{N},\,u\in U^{(N)}\right)\,\,\,\,\, \text{and}\,\,\,\,\, Y=\left(Y^{(N)}(u)|N\in\mathbb{N},\,u\in U^{(N)}\right)\] are families of non negative random variables indexed by $N$, and possibly some parameter $u$, then we say that $X$ is stochastically dominated by $Y$, if for all $\epsilon,D>0$ we have \[\sup_{u\in U^{(N)}}\mathbb{P}\left(X^{(N)}(u)>N^\epsilon Y^{(N)}(u)\right)\le N^{-D}\] for large enough $N\ge N_0(\epsilon,D)$. In this case we use the notation $X\prec Y$. Moreover, if we have $|X|\prec Y$, we also write $X=\mathcal{O}_\prec(Y)$.
\end{definition}

We will say that a sequence of events $A=A^{(N)}$ holds with overwhelming probability if $\mathbb{P}\left(A^{(N)}\right)\ge 1-N^{-D}$ for any $D>0$ and $N\ge N_0(D)$. In particular, under the conditions \eqref{var} and \eqref{mom}, we have $X_{i\mu}\prec (MN)^{\frac{1}{4}}$ uniformly in $i,\mu$ and that $\max_k \lambda_k\le \gamma_+ +\delta_*/2$, $\min_k \lambda_k \ge \max\{0,\gamma_- -\delta_*/2\}$ with overwhelming probability (see Theorem 2.10, Lemma 4.11 in \cite{bloemendal2014isotropic}).

We define the almost analytic extension of $f\in H_0^2([\gamma_--\delta_*,\gamma_++\delta_*])$ by
\begin{equation}
\label{aae}
f_\mathbb{C}(x+i\eta):=\left(f(x)+i\eta f'(x)\right)\tilde{\chi}(\eta),
\end{equation}
where $\tilde{\chi}:\mathbb{R}\to\mathbb{R}$ is a smooth cut-off function which is constant 1 in $[-5,5]$ and constant 0 outside $[-10,10]$. By this definition it follows that $f_\mathbb{C}$ is bounded and compactly supported, i.e. $f_\mathbb{C}(x+i\eta)=0$ for $x\in [\gamma_--\delta_*,\gamma_++\delta_*]^c$. Furthermore for small $\eta$ we have that \begin{equation}
\label{fsmall}
 \partial_{\overline{z}}f_\mathbb{C}(x+i\eta)=\mathcal{O}(\eta)\,\,\, \text{and}\,\,\, \partial_\eta\partial_{\overline{z}}f_\mathbb{C}(x+i\eta)=\mathcal{O}(1).
\end{equation}


We use the following representation of $f_N$ from \cite{erdHos2017fluctuations}:\begin{equation}
\label{formula}
f_N=\frac{2}{\pi}\Re \int_\mathbb{R}\int_{\mathbb{R}_+}\partial_{\overline{z}}f_\mathbb{C}(x+i\eta)\Delta_N(x+i\eta)\,dxd\eta.
\end{equation}

We first exclude a critical area very close to the real line in the integral in \eqref{formula}. From the resolvent identities 
$\left|\eta  \left\langle {\bf x},G^2{\bf x}\right\rangle\right|\le \Im\left\langle {\bf x},G {\bf x}\right\rangle.$
Then, we have that \[\left|\eta z\left\langle {\bf x},G^2{\bf x}\right\rangle+\eta\left\langle {\bf x},G{\bf x}\right\rangle+\eta\right|\le 2 \left|z+z\left\langle {\bf x},G{\bf x}\right\rangle\right|.\] 
Hence, we conclude that \begin{equation}
\label{bound}
\left|\eta\Delta_N(x+i\eta)\right|\le 2.
\end{equation}

To study $f_N$ we restrict our integration to the domain $\Im z\in [\eta_0,10]$, with \\ $\eta_0:=N^{-\frac{2}{3}}$. Thanks to \eqref{fsmall} and \eqref{bound}, we find that \[f_N=\frac{2}{\pi}\Re \int_\mathbb{R}\int_{\eta_0}^{10}\partial_{\overline{z}}f_\mathbb{C}(x+i\eta)\Delta_N(x+i\eta)\,dxd\eta+\mathcal{O}_\prec\left(\eta_0\right).\] Then, for $\Im z=\eta\ge \eta_0$ we claim that the leading order term of $\Delta_N(z)$ is given by \begin{equation}
\label{leading}
\widehat{\Delta}_N(z):=\frac{1+\frac{1}{N\sqrt{\phi}}\text{Tr}G(z)+\tilde{z}\frac{1}{N}\text{Tr}G^2(z)}{-z-\tilde{z}\frac{1}{N}\text{Tr}G(z)},
\end{equation} with the notation $\tilde{z}:=\phi^{-\frac{1}{2}}z$ for brevity. Note that \eqref{leading} is related to \eqref{delta12} by taking expectation with respect to ${\bf x}$ in the numerator and denominator separately. 

We split the analysis of $f_N$ into two parts: the leading order term 
\begin{equation}
\label{leadingt}
\widehat{\Omega}_f:=\frac{2}{\pi}\Re \int_\mathbb{R}\int_{\eta_0}^{10}\partial_{\overline{z}}f_\mathbb{C}(x+i\eta)\widehat{\Delta}_N(x+i\eta)\, d\eta dx
\end{equation}
and the fluctuation term \begin{equation}
\label{fluctuation}
F_N:=\frac{2}{\pi}\Re \int_\mathbb{R}\int_{\eta_0}^{10}\partial_{\overline{z}}f_\mathbb{C}(x+i\eta)\left(\Delta_N(x+i\eta)-\widehat{\Delta}_N(x+i\eta)\right)\, d\eta dx.
\end{equation} In this way we have that \[f_N=\widehat{\Omega}_f+F_N+\mathcal{O}_\prec\left(N^{-\frac{2}{3}}\right).\]
In the following two sections, we will show that $\widehat{\Omega}_f=\Omega_f+\mathcal{O}_\prec(N^{-\frac{2}{3}})$ and $\mathbb{E}(F_N^2)=\frac{1}{N}V_f+\mathcal{O}_\prec(N^{-\frac{7}{6}})$, with some $N$-independent constant $V_f$, which will prove Theorem \ref{maintheo} in the sense of mean and variance.

\subsection{Leading term: calculation of  the mean.}
The main tool we will use is the local law for the Marchenko-Pastur distribution in its averaged and entry-wise from. These results have first been proven in \cite{bloemendal2014isotropic} (see Theorem 2.4 and Theorem 2.5) uniformly for each $z\in\bf{S}$, where  \[
\textbf{S}\equiv \textbf{S}(\omega,\eta_0):=\left\{z=x+i\eta\in\mathbb{C}:\kappa_x\le \omega^{-1},\, \eta_0\le \eta\le \omega^{-1},\, |z|\ge \omega\right\},
\] with some $\omega\in(0,1)$ fixed and $\kappa_x:=\min\{|\gamma_+-x|,\,|\gamma_--x|\}$. In our proof, instead, we rely on local laws which hold true for each $z\in\mathbb{H}$, hence, combining the results in \cite{bloemendal2014isotropic} with Theorem 2.7 and Theorem 2.9 respectively for $\phi=1$ and $d_*\le|\phi-1|\le \hat{d}$ in \cite{alt2017local}, we get the Marchenko-Pastur local law in the averaged form \begin{equation}
\label{average}
\begin{split}
&m_R(z):=\frac{1}{N}\text{Tr}R(z)=m_\phi(z)+\mathcal{O}_\prec \left(\frac{1}{N\eta}\right),\\
&m_G(z):=\frac{1}{M}\text{Tr}G(z)=m_{\phi^{-1}}(z)+\mathcal{O}_\prec\left(\frac{1}{N\eta}\right),\end{split}
\end{equation} and its entry-wise form \begin{equation}
\label{entrywise}
\left|R_{\mu\nu}(z)-\delta_{\mu\nu}m_\phi(z)\right|\prec \frac{1}{\sqrt{N\eta |z|}}, \,\,\,\left|G_{ij}(z)-\delta_{ij}m_{\phi^{-1}}(z)\right|\prec\frac{1}{\sqrt{N\eta |z|}}
\end{equation} uniformly for each $z\in\mathbb{H}$.

\begin{remark}
Notice that in \eqref{average} and \eqref{entrywise} the error term from \cite{alt2017local} is smaller in some particular cases, but we will not need these optimal bounds and we write local laws in a unified form which hold true for both the cases $\phi=1$ and $d_*\le|\phi-1|\le \hat{d}$.
\end{remark}

By \eqref{tracem}, we have that 
\begin{equation}
\label{GRtr}
\tilde{z}m_G(z)=\tilde{z}m_R(z)-\phi^\frac{1}{2}+\phi^{-\frac{1}{2}}+\frac{\phi^\frac{1}{2}}{N}.
\end{equation}

Hence, using the equality above, \eqref{st1} and \eqref{average}, we write \eqref{leading} as follows
\begin{equation}
\label{good}
\begin{split}\widehat{\Delta}_N(z)&=\frac{1+\frac{1}{N\sqrt{\phi}}\text{Tr}G(z)+\tilde{z}\frac{1}{N}\text{Tr}G^2(z)}{-z+\phi^\frac{1}{2}-\phi^{-\frac{1}{2}}-\tilde{z}m_R(z)}\\
&=m_\phi(z)\left(1+\frac{1}{N\sqrt{\phi}}\text{Tr}G(z)+\tilde{z}\frac{1}{N}\text{Tr}G^2(z)\right)+\mathcal{O}_\prec\left(\frac{1}{N\eta}\right).\end{split}
\end{equation}

Then, thanks to \eqref{good} and \eqref{fsmall}, we obtain \[ \widehat{\Omega}_f=\frac{2}{\pi}\Re \int_\mathbb{R}\int_{\eta_0}^{10}\partial_{\overline{z}}f_\mathbb{C}(z)m_\phi(z)\left(1+\frac{1}{N\sqrt{\phi}}\text{Tr}G(z)+\frac{\tilde{z}}{N}\text{Tr}G^2(z)\right)d\eta dx+\mathcal{O}_\prec\left(\frac{1}{N}\right),\] where from now on we will use the notation $z=x+i\eta$ and $z_0=x+i\eta_0$. Furthermore, we notice that, using \eqref{GRtr} and the identity $\partial_z \text{Tr}G(z)=\text{Tr}G^2(z)$, we get \[\begin{split}1+\frac{1}{N\sqrt{\phi}}\text{Tr}G(z)+\tilde{z}\frac{1}{N}\text{Tr}G^2(z)&=\partial_\eta\left(\eta-i\phi^{-\frac{1}{2}}(x+i\eta)\frac{1}{N}\text{Tr}G(x+i\eta)\right) \\
&=\partial_\eta\left(\eta-i\phi^{-\frac{1}{2}}(x+i\eta)m_R(x+i\eta)\right).\end{split}\] 
Hence, integrating by parts twice in $\eta$, using that the upper limit of the $\eta-$integration is zero since $\partial_{\overline{z}}f_\mathbb{C}(x+10i)=0$ by the definition of $\tilde{\chi}$, we have 
\begin{align*}
\widehat{\Omega}_f&=-\frac{2}{\pi}\Re \int_\mathbb{R}\partial_{\overline{z}}f_\mathbb{C}(z_0)m_\phi(z_0)(\eta_0-i\tilde{z}_0 m_R(z_0))\, dx\\
&\quad-\frac{2}{\pi}\Re \int_\mathbb{R}\int_{\eta_0}^{10}\partial_\eta\big( \partial_{\overline{z}}f_\mathbb{C}(z)m_\phi(z)\big)(\eta-i\tilde{z}m_R(z))\, d\eta dx+\mathcal{O}_\prec\left( N^{-1}\right) \\
&=-\frac{2}{\pi}\Re \int_\mathbb{R}\int_{\eta_0}^{10}\partial_\eta \big(\partial_{\overline{z}}f_\mathbb{C}(z)m_\phi(z)\big)(\eta-i\tilde{z}m_\phi(z))\, d\eta dx \numberthis\label{estimatelt}\\
&\quad-\frac{2}{\pi}\Re \int_\mathbb{R}\int_{\eta_0}^{10}\partial_\eta \big(\partial_{\overline{z}}f_\mathbb{C}(z)m_\phi(z)\big)(-i\tilde{z}m_R(z)+i\tilde{z}m_\phi(z))\, d\eta dx+\mathcal{O}_\prec\left(\eta_0\right) \\
&=\frac{2}{\pi}\Re \int_\mathbb{R}\int_{\eta_0}^{10}\partial_{\overline{z}}f_\mathbb{C}(z)m_\phi(z)(1+(\tilde{z}m_\phi(z))')d\eta dx+\mathcal{O}_\prec\left(\eta_0\right)+\mathcal{O}_\prec\left(\frac{|\log \eta_0|}{N}\right),
\end{align*}
where we used that $\partial_{\overline{z}}f_\mathbb{C}(x+i\eta)$ scales like $\eta$ near the real axis by \eqref{fsmall}, the local law from \eqref{average} and that $|z\phi^{-\frac{1}{2}}\partial_\eta\big(\partial_{\overline{z}}f_\mathbb{C}m_\phi(z)\big)|\le C$ from the bounds \eqref{goodbound} and \eqref{boundderivative}. In the last step we also used that $-i\partial_\eta h(z)=\partial_z h(z)$ for any analytic function $h$.

In summary, by \eqref{estimatelt}, we conclude that
\begin{equation}
\label{newomegaf}
\widehat{\Omega}_f=\frac{2}{\pi}\Re \int_\mathbb{R}\int_{\eta_0}^{10}\partial_{\overline{z}}f_\mathbb{C}(z)p_\phi(z)\,d\eta dx+\mathcal{O}_\prec(\eta_0),
\end{equation} where for brevity we introduced
\begin{equation}
\label{newp}
p_\phi(z):=m_\phi(z)[1+(\tilde{z}m_\phi(z))'],\qquad z\in\mathbb{H}.
\end{equation}

For the main term we need the following lemma (see Lemma 3.4  in \cite{erdHos2017fluctuations}).

\begin{lemma}
\label{sto}
Let $\varphi,\psi:[a,b]\times [0,10i]\to\mathbb{C}$ be functions such that $\partial_{\overline{z}}\psi(z)=0$, $\varphi,\psi\in H^1$ and $\varphi$ vanishes at the left, right and top of the boundary of the integration region. Then for any $\tilde{\eta}\in [0,10]$, we have \[\int_a^b\int_{\tilde{\eta}}^{10}(\partial_{\overline{z}}\varphi(z))\psi(z)\, d\eta dx=\frac{1}{2i}\int_a^b\varphi(x+i\tilde{\eta})\psi(x+i\tilde{\eta})\, dx.\]
\end{lemma}

In order to compute the leading term defined in \eqref{leadingt} we extend the integral in \eqref{newomegaf} to the real axis.
 For this purpose we introduce a tiny auxiliary scale $\eta_1$, say $\eta_1:=N^{-10}$. We recall that $f_\mathbb{C}$
  is supported in $[a,b]\times [-10,10]$, with $a=\gamma_--\delta_*$ and $b=\gamma_++\delta_*$, where $\gamma_-,\gamma_+$ are the spectral edges, and $\kappa_x=\min\{|x-\gamma_-|,|x-\gamma_+|\}$.

Since by \eqref{goodbound}, \eqref{boundderivative} and \eqref{fsmall}, we have that 
\[
\left|\frac{2}{\pi}\Re \int_\mathbb{R}\int_{\eta_1}^{\eta_0}\partial_{\overline{z}}f_\mathbb{C}(z)p_\phi(z)\,d\eta dx\right|\lesssim \int_a^b\int_{\eta_1}^{\eta_0}\left(\frac{\eta}{|z|}+\frac{\eta}{|z|^\frac{1}{2}\sqrt{\kappa_x+\eta}}\right)\,d\eta dx\lesssim\eta_0^\frac{3}{2},
\] we conclude that 
\begin{equation}
\label{newlt}
\widehat{\Omega}_f=\frac{2}{\pi}\Re \int_\mathbb{R}\int_{\eta_1}^{10}\partial_{\overline{z}}f_\mathbb{C}(z)p_\phi(z)\,d\eta dx+\mathcal{O}_\prec\left(\eta_0\right).
\end{equation} 

Next, applying Lemma \ref{sto} to the integral in the r.h.s. of \eqref{newlt}, we conclude
\begin{equation}
\widehat{\Omega}_f=\frac{1}{\pi}\Im \int_\mathbb{R}f_\mathbb{C}(x+i \eta_1)p_\phi(x+i\eta_1)\, dx+\mathcal{O}_\prec\left(\eta_0\right).
\end{equation} 

By \eqref{aae} and \eqref{fsmall}, using the bounds \eqref{goodbound}--\eqref{boundderivative}, it easily follows that 
\begin{equation}
\label{cona}
\widehat{\Omega}_f=\frac{1}{\pi}\int_\mathbb{R}f(x)\Im p_\phi(x+i\eta_1)\, dx+\mathcal{O}_\prec\left(\eta_0\right).
\end{equation}

We notice that \begin{equation}
\label{wigner1}
w_\phi(z):=\sqrt{\phi}(1+zm_{\phi^{-1}}(z))=\frac{\phi^\frac{1}{2}+\phi^{-\frac{1}{2}}-z+i\sqrt{(z-\gamma_-)(\gamma_+-z)}}{2}
\end{equation} is the Stieltjes transform of the Wigner semicircle law centered at $\phi^\frac{1}{2}+\phi^{-\frac{1}{2}}$. Hence, $w_\phi$ is also characterized as the unique solution of \begin{equation}
\label{stw}
w_\phi(z)+\frac{1}{z-\phi^\frac{1}{2}-\phi^{-\frac{1}{2}}+w_\phi(z)}=0, \,\,\,\,\, \Im w_\phi>0.
\end{equation} Notice that $w_\phi(z)=w_{\phi^{-1}}(z)$ and that, using the self consistent equation \eqref{st1} and the relation between $m_\phi$ and $m_{\phi^{-1}}$ in \eqref{change1}, we have \begin{equation}
\label{wwin}
w_\phi(z)=-zm_\phi(z)m_{\phi^{-1}}(z).
\end{equation}

We now distinguish the cases $\phi=1$ and $|\phi-1|\ge d_*$, since for $\phi=1$ the integral in \eqref{cona} has an additional singularity in zero which we have to take into account. 

We start with the case $|\phi-1|\ge d_*$. In this case $\gamma_-\ge \tau(d_*)$, for some $\tau(d_*)>0$. By equations \eqref{st1} and \eqref{change1}, expressing $(\tilde{z}m_\phi)'=w_\phi'$ from differentiating the self consistent equation for $w_\phi$ in \eqref{stw}, it follows that 
\begin{equation}
\label{derwder}
w_\phi(z)'=\frac{w_\phi^2(z)}{1-w_\phi^2(z)},
\end{equation} and so we may write $p_\phi$ from \eqref{newp} as 
 \begin{equation}
\label{newpp}
p_\phi(z)=\frac{m_\phi(z)}{1-w_\phi^2(z)}.
\end{equation}

Furthermore, by Lemma 3.6 of \cite{bloemendal2016principal} we have that
\begin{equation}
\label{derw}
|1-w_\phi^2(z)|\asymp \sqrt{\kappa_x+\eta},\qquad c\le |w_\phi(z)|\le 1,
\end{equation} 
with some $\phi$-independent constant $c>0$, for any $z=x+i\eta$ such that $|z-\sqrt{\phi}|\le 10$.

To evaluate $\widehat{\Omega}_f$ in \eqref{cona}, we first remove the $\eta_1$ in the argument of $p_\phi$. We proceed writing $p_\phi(x+i\eta_1)-p_\phi(x)$ as follows 
\begin{align*}
p_\phi(x+i\eta_1)-p_\phi(x)&=\frac{1}{1-w_\phi(x+i\eta_1)^2}\int_0^{\eta_1} m_\phi(x+i\eta)'\, d\eta \numberthis \label{longform}\\
&\quad +\frac{m_\phi(x) (w_\phi(x+i\eta_1)+w_\phi(x))}{(1-w_\phi(x+i\eta_1)^2)(1-w_\phi(x)^2)}\int_0^{\eta_1} w_\phi(x+i\eta)'\, d\eta.
\end{align*}

Then, by \eqref{goodbound}--\eqref{boundderivative} and \eqref{derwder}--\eqref{longform}, simple estimates give that 
\begin{align*}
\Big|p_\phi(x+i\eta_1)-p_\phi(x)\Big|&\lesssim \frac{\eta_1^{1/4}}{|x|^{1/2}\kappa_x^{1/4}\sqrt{\kappa_x+\eta_1}}+\frac{\sqrt{\eta_1}}{|x|^{1/2}\sqrt{\kappa_x(\kappa_x+\eta_1)}}\lesssim \frac{\eta_1^{1/4}}{|x|^{1/2}\kappa_x^{3/4}},
\end{align*} for any $x\in \mathbb{R}$. Hence, if $|\phi-1|\ge d_*$, integrating over $x$, we conclude that 
\begin{equation}
\label{firstca}
\left|\frac{1}{\pi}\int_a^bf(x)\Im [p_\phi(x+i\eta_1)-p_\phi(x)]\, dx\right|\lesssim \eta_1^{1/4}.
\end{equation}
In particular, this implies that $\Im p_\phi(z)$
is of order $\eta_1^{1/4}$ outside the interval $[\gamma_-,\gamma_+]$, since  $\Im p_\phi(x)=0$ for $x\notin [\gamma_-,\gamma_+]$.
Moreover, \eqref{newp},  \eqref{cona} and \eqref{firstca} imply that 
\[ 
\begin{split}
\widehat{\Omega}_f=\frac{1}{\pi}\int_{\gamma_-}^{\gamma_+}f(x)\Im \left[m_\phi(x)(1+(x\phi^{-\frac{1}{2}}m_\phi(x))')\right]\,dx+\mathcal{O}_\prec\left(N^{-\frac{2}{3}}\right) \\
=\int_{\gamma_-}^{\gamma_+} f(x)\frac{\sqrt{\phi}}{4\pi^2 x\rho_\phi(x)}\left(1+\frac{\sqrt{\phi}-\frac{1}{\sqrt{\phi}}}{x}\right)\,dx+\mathcal{O}_\prec\left(N^{-\frac{2}{3}}\right),
\end{split}
\] 
concluding the estimate for the leading term of $\mathbb{E}f_N$ when $|\phi-1|\ge d_*$. 

Now we consider the case $\phi=1$, when $\gamma_-=0$ and $\gamma_+=4$. In this case, the computation of the integral \eqref{cona} is a bit more delicate since the singularities around $x\approx 0$ 
and $\kappa_x\approx 0$ overlap. For brevity, in the rest  of this section we use the notation $m=m(z):=m_{\phi=1}(z)$ and $w=w(z):=w_{\phi=1}(z)$ for any $z\in \mathbb{H}$. Expressing $m'$ from differentiating the self consistent equation \eqref{st1}, using \eqref{st1} repeatedly and the relation \eqref{wwin}, a simple calculation gives that 
\begin{equation}
\label{newformm}
p=m(1+(zm)')=-\frac{1}{z}\cdot \frac{1}{1-zm^2}=-\frac{1}{z}\cdot\frac{1}{1+w},
\end{equation} with $p=p(z):=p_{\phi=1}(z)$.
We also define 
\begin{equation}
\label{newq}
q(z):=\frac{1}{1+w(z)}, \qquad z\in \mathbb{H}.
\end{equation}

As a consequence of \eqref{newformm}-\eqref{newq}, it follows that 
\begin{align*}
&\left|\int_a^bf(x)\Im [p(x+i\eta_1)-p(x)]\, dx-\frac{\pi f(0)}{2}\right|\\
&\qquad\qquad\le \left|\int_a^bf(x)\Big[\frac{x}{x^2+\eta_1^2}\Im q(x+i\eta_1)-\frac{1}{x}\Im q(x)\Big]\, dx\right| \numberthis \label{ima}\\
&\qquad\qquad\quad +\left|\int_a^bf(x)\frac{\eta_1}{x^2+\eta_1^2} \Re q(x+i\eta_1)\, dx-\frac{\pi f(0)}{2}\right|. \numberthis \label{rea}
\end{align*}
We start estimating \eqref{rea}. Using explicit computations, by the expression in \eqref{wigner1} for $\phi=1$, we conclude that 
\begin{equation}
\label{est1}
\eqref{rea}\le \left|\frac{1}{2}\int_a^b f(x)\frac{\eta_1}{x^2+\eta_1^2}\, dx-\frac{\pi f(0)}{2}\right|+\mathcal{O}(\sqrt{\eta_1})\lesssim \sqrt{\eta_1}.
\end{equation}
Furthermore, since
\[
|1+w(z)|=|1-zm(z)^2|\asymp \frac{\sqrt{\kappa_x+\eta}}{|z|^\frac{1}{2}},
\]
by \eqref{finaleq}, using \eqref{wigner1} and the definition of $q$ in \eqref{newq}, it also follows that the integrand in \eqref{ima} is bounded by 
\begin{equation}
\label{est2}
\frac{f(x)|x|^{3/2}\sqrt{\eta_1}}{(x^2+\eta_1^2)^\frac{3}{4}\sqrt{\kappa_x(\kappa_x+\eta_1)}}+\frac{f(x)\eta_1^2}{|x|^{1/2}|4-x|^{1/2}(x^2+\eta_1^2)},
\end{equation} for any $x\in \mathbb{R}$.
Then, combining \eqref{est1} with the integral of \eqref{est2}, we conclude
\begin{equation}
\label{secondca}
\left|\frac{1}{\pi}\int_a^b f(x)\Im [p(x+i\eta_1)-p(x)]\, dx-\frac{f(0)}{2}\right|\lesssim \eta_1^{1/4}.
\end{equation}
Similarly to the case $|\phi-1|\ge d_*$, this bound implies that $\Im p(x+i\eta_1)$ is of order $\eta_1^{1/4}$ outside $[0,4]$. Hence, the above inequality implies that
\[
\widehat{\Omega}_f=\frac{1}{\pi}\int_0^4 \frac{f(x)}{4\pi^2 x\rho_1(x)}\,dx+\frac{f(0)}{2}+\mathcal{O}_\prec\left(N^{-\frac{2}{3}}\right),
\] concluding the computation of $\Omega_f$, the leading term of $\mathbb{E} f_N$ in Theorem~\ref{maintheo}.

\subsection{Fluctuation term}
We write the difference $\Delta_N(z)-\widehat{\Delta}_N(z)$ in a more convenient form to study the integral in \eqref{fluctuation}. The key point is to express it as a derivative (up to an error) to prepare it for an integration by parts. Let $\hat{z}$ be defined as $\hat{z}:=z\phi^\frac{1}{2}$.

\begin{lemma}
For any $\eta>\eta_0$ we have that \begin{equation}
\label{der}
\Delta_N(z)-\widehat{\Delta}_N(z)=\partial_z\frac{z\left\langle {\bf x},G{\bf x}\right\rangle-\hat{z}m_G(z)}{-z-\hat{z}m_G(z)}+\mathcal{O}_\prec\left(\frac{1}{N\eta^2}\right).
\end{equation}
\proof
\normalfont
This lemma, using \eqref{mom}, relies on the following large deviation bound (see, e.g. Lemma 3.1 in \cite{bloemendal2014isotropic}) \begin{equation}
\label{large}
\left\langle {\bf x},G{\bf x}\right\rangle=\frac{1}{\sqrt{MN}}\text{Tr}G+\mathcal{O}_\prec\left(\sqrt{(MN)^{-1}\text{Tr}|G|^2}\right),
\end{equation} and a similar formula for $\left\langle {\bf x},G^2{\bf x}\right\rangle$.

In the following part of the proof, in order to abbreviate our notation, we use $G:=G(z)$, $m_G:=m_G(z)$. Using \eqref{delta12} and \eqref{leading}, we have \begin{equation}
\label{alpha}
\begin{split}\Delta_N(z)-\widehat{\Delta}_N(z)&=\frac{\left(z\left\langle {\bf x},G^2{\bf x}\right\rangle+\left\langle {\bf x},G{\bf x}\right\rangle+1\right)(-z-\hat{z}m_G)}{\left(-z-z\left\langle {\bf x},G {\bf x}\right\rangle\right)(-z-\hat{z}m_G)}\\
&\quad+\frac{(-1-\phi^\frac{1}{2}m_G-\hat{z}m_G')\left(-z-z\left\langle {\bf x},G {\bf x}\right\rangle\right)}{\left(-z-z\left\langle {\bf x},G {\bf x}\right\rangle\right)(-z-\hat{z}m_G)}.\end{split}
\end{equation}

Now we claim that \[\Delta_N(z)-\widehat{\Delta}_N(z)=\partial_z\frac{z\left\langle {\bf x},G{\bf x}\right\rangle-\hat{z}m_G}{-z-\hat{z}m_G}+\mathcal{E},\] with an error term $\mathcal{E}$ we will determine along the proof. We start with \begin{align*}
\partial_z\frac{z\left\langle {\bf x},G{\bf x}\right\rangle-\hat{z}m_G}{-z-\hat{z}m_G}&=\frac{(-z-\hat{z}m_G)\left(\left\langle {\bf x},G{\bf x}\right\rangle+z\left\langle {\bf x},G^2{\bf x}\right\rangle-\phi^\frac{1}{2}m_G-\hat{z}m_G'\right)}{(-z-\hat{z}m_G)^2} \\
&\quad-\frac{(-1-\phi^\frac{1}{2}m_G-\hat{z}m_G')\left(z\left\langle {\bf x},G{\bf x}\right\rangle-\hat{z}m_G\right)}{(-z-\hat{z}m_G)^2}.\numberthis\label{beta}
\end{align*} Using $m_G(z)=\frac{1}{M}\text{Tr}G(z)$ and $m_G'(z)=\frac{1}{M}\text{Tr}G^2(z)$ we write the r.h.s. of \eqref{alpha} as \begin{equation}\begin{split}
\label{gamma}
\Delta_N(z)-\widehat{\Delta}_N(z)&=\frac{\left\langle {\bf x},G{\bf x}\right\rangle+z\left\langle {\bf x},G^2{\bf x}\right\rangle-\phi^\frac{1}{2}m_G-\hat{z}m_G'}{(-z-\hat{z}m_G)-\left(z\left\langle {\bf x},G{\bf x}\right\rangle-\hat{z}m_G\right)} \\
&\quad-\frac{(-1-\phi^\frac{1}{2}m_G-\hat{z}m_G')\left(z\left\langle {\bf x},G{\bf x}\right\rangle-\hat{z}m_G\right)}{(-z-\hat{z}m_G)^2-(-z-\hat{z}m_G)\left(z\left\langle {\bf x},G{\bf x}\right\rangle-\hat{z}m_G\right)}.\end{split}
\end{equation}
By \eqref{average}, \eqref{large} and the bound in \eqref{goodbound} it follows that \begin{equation}
\label{average1}
z\left\langle {\bf x},G{\bf x}\right\rangle-\hat{z}m_G(z)\prec \frac{|z|}{\sqrt{MN}}\sqrt{\text{Tr}|G(z)|^2}\le \frac{|z|}{\sqrt{MN}}\sqrt{\frac{1}{\eta}\Im \text{Tr}G(z)}\prec\frac{|z|^\frac{3}{4}}{\sqrt{N\eta}}
\end{equation} and also \begin{equation}
\label{average2}
z\left\langle {\bf x},G^2{\bf x}\right\rangle-\hat{z}m_G'(z)\prec \frac{|z|}{\sqrt{MN}}\sqrt{\text{Tr}|G|^4}\le \frac{|z|}{\sqrt{MN}\eta}\sqrt{\text{Tr}|G(z)|^2}\prec\frac{|z|^\frac{1}{2}}{\sqrt{N\eta^3}}.
\end{equation}  

Note that the leading term in the denominators in \eqref{gamma} is separated away from zero since $-z-\hat{z}m_{\phi^{-1}}(z)=[m_{\phi^{-1}}(z)]^{-1}$, by \eqref{st1} and \eqref{change1}. Thus these denominators are stable under small perturbations. Hence, replacing $z\left\langle {\bf x},G{\bf x}\right\rangle$ in the denominator with $\hat{z}m_G(z)+\mathcal{O}_\prec\left(\frac{1}{\sqrt{N\eta}}\right)$ and comparing \eqref{beta} and \eqref{gamma}, we conclude that \[\Delta_N(z)-\widehat{\Delta}_N(z)=\partial_z\frac{z\left\langle {\bf x},G{\bf x}\right\rangle-\hat{z}m_G}{-z-\hat{z}m_G}+\mathcal{O}_\prec\left(\frac{1}{N\eta^2}\right).\] In estimating various error terms along the proof we used that $zm_G(z)=\mathcal{O}_\prec(1)$ (by \eqref{average} and \eqref{goodbound}) and that $zm_G'(z)=\mathcal{O}_\prec\left(\eta^{-1}\right)$ by \eqref{ward} and \eqref{goodbound}.
\endproof
\end{lemma}

Next, we use \eqref{der} to estimate the fluctuation term $F_N$ as defined in \eqref{fluctuation} via an integration by parts \[\begin{split}F_N&= -\frac{2}{\pi}\Re\int_\mathbb{R}\partial_{\overline{z}}f_\mathbb{C}(z_0)i\frac{z_0\left\langle {\bf x},G(z_0){\bf x}\right\rangle-\hat{z}_0m_G(z_0)}{-z_0-\hat{z}_0m_G(z_0)}\, dx \\
&\quad+\frac{2}{\pi}\Re \int_\mathbb{R}\int_{\eta_0}^{10}\partial_\eta\partial_{\overline{z}}f_\mathbb{C}(z)i\frac{z\left\langle {\bf x},G(z){\bf x}\right\rangle-\hat{z}m_G(z)}{-z-\hat{z}m_G(z)}\,d\eta dx+\mathcal{O}_\prec\left(\frac{|\log \eta_0|}{N}\right),\end{split}\] with $\hat{z}_0:=\phi^\frac{1}{2}z_0$. Then, we continue with the estimate \[\frac{z\left\langle {\bf x},G(z){\bf x}\right\rangle-\hat{z}m_G(z)}{-z-\hat{z}m_G(z)}=m_\phi(z)\left(z\left\langle {\bf x},G(z){\bf x}\right\rangle-\hat{z}m_G(z)\right)+\mathcal{O}_\prec\left(\frac{1}{(N\eta)^\frac{3}{2}}\right)\] from \eqref{average}, \eqref{trace}, \eqref{st1} and \eqref{average1} to find that \begin{align*} 
F_N&= -\frac{2}{\pi}\Re \int_\mathbb{R} m_\phi(z_0)\partial_{\overline{z}}f_\mathbb{C}(z_0)i\left(z_0\left\langle {\bf x},G(z_0){\bf x}\right\rangle-\hat{z}_0m_G(z_0)\right)\, dx \numberthis\\
&\quad+\frac{2}{\pi}\Re \int_\mathbb{R}\int_{\eta_0}^{10} m_\phi(z)\partial_\eta\partial_{\overline{z}}f_\mathbb{C}(z)i\left(z\left\langle {\bf x},G(z){\bf x}\right\rangle-\hat{z}m_G(z)\right)\,d\eta dx+\mathcal{O}_\prec\left(N^{-\frac{2}{3}}\right) \\
&=-\frac{2}{\pi}\Im \int_\mathbb{R}\int_{\eta_0}^{10} m_\phi(z)\partial_\eta\partial_{\overline{z}}f_\mathbb{C}(z)\left(z\left\langle {\bf x},G(z){\bf x}\right\rangle-\hat{z}m_G(z)\right)\,d\eta dx+\mathcal{O}_\prec\left(N^{-\frac{2}{3}}\right) ,\end{align*} where in the last step we used that by \eqref{fsmall} and \eqref{average1} it follows \[\left|\partial_{\overline{z}}f_\mathbb{C}(z_0)i\left(z_0\left\langle {\bf x},G(z_0){\bf x}\right\rangle-\hat{z}_0m_G(z_0)\right)\right|\prec \sqrt{\frac{\eta_0}{N}}\le N^{-\frac{2}{3}}.\]

The leading order expression for $F_N$ has zero mean, hence we can start computing the variance $\text{Var}(F_N)=\mathbb{E}F_N^2+\mathcal{O}_\prec\left(N^{-\frac{4}{3}}\right)$ as 
\[
\mathbb{E}F_N^2=\mathbb{E}\left(\frac{2}{\pi}\Im \int_\mathbb{R}\int_{\eta_0}^{10} m_\phi(z)\partial_\eta\partial_{\overline{z}}f_\mathbb{C}(z)\left(z\left\langle {\bf x},G(z){\bf x}\right\rangle-\hat{z}m_G(z)\right)d\eta dx\right)^2+\mathcal{O}_\prec\left(N^{-\frac{7}{6}}\right).  
\]

When we use the expectation $\mathbb{E}$ we frequently use the property that if $X$ and $Y$ are random variables with $X=\mathcal{O}_\prec(Y)$, $Y\ge 0$ and $|X|\le N^C$ for some constant $C$, then $\mathbb{E}|X|\prec\mathbb{E}Y$, or, equivalently, $\mathbb{E}|X|\le N^\epsilon\mathbb{E}Y$ for any $\epsilon>0$ and $N\ge N_0(\epsilon)$. 
To compute the leading term $F_N'$ in $\mathbb{E}F_N^2$ we introduce the short-hand notations \begin{equation}
\label{short}
g(z):=\frac{2}{\pi}z m_\phi(z) \partial_\eta\partial_{\overline{z}}f_\mathbb{C}(z),\,\,\, A(z):=\sqrt{N}\left(\left\langle {\bf x},G(z){\bf x}\right\rangle- \phi^\frac{1}{2}m_G(z)\right)
\end{equation} to write 
\[ F_N':=\frac{1}{N}\mathbb{E}\left(\Im\int_\mathbb{R}\int_{\eta_0}^{10} g(z)A(z)\,d\eta dx\right)^2.\] 
We will often use the following identity for any $z,w\in\mathbb{C}$: 
\begin{equation}
\label{complex}
(\Im z)(\Im w)=\frac{1}{2}\Re(\overline{z}w-zw).
\end{equation} 
Thanks to \eqref{complex} we write 
\begin{equation}
\label{var56}
F_N'=\frac{1}{2N}\Re \iint_\mathbb{R}\iint_{\eta_0}^{10}\big[ g(z)g(\overline{z}')\mathbb{E}\left(A(z)A(\overline{z}')\right)-g(z)g(z')\mathbb{E}\left(A(z)A(z')\right)\big]\,d\eta d\eta' dxdx',
\end{equation} where we used that $\overline{X(z)}=X(\overline{z})$ and $\overline{g(z)}=g(\overline{z})$. In the following we use the short notation $G=G(z), \,\, G'=G(z')$.

To study the expectation of $A(z)A(z')$, we consider \[\begin{split}A(z)A(z')=&N\left(\sum_{i,j=1,\,i\ne j}^M\overline{{\bf x}_i}G_{ij}{\bf x}_j+\sum_{i=1}^M \left( |{\bf x}_i|^2-\frac{1}{\sqrt{MN}}\right)G_{ii}\right)\\
&\times\left(\sum_{l,k=1,\, l\ne k}^M\overline{{\bf x}_l}G_{lk}'{\bf x}_k+\sum_{l=1}^M \left( |{\bf x}_l|^2-\frac{1}{\sqrt{MN}}\right)G_{ll}'\right).\end{split}\] The conditional expectation $\mathbb{E}_1=\mathbb{E}(\cdot|X)$ conditioned on the matrix $X$ gives \begin{align*}
\mathbb{E}_1(A(z)A(z'))&=\frac{1}{\phi N}\sum_{i,j=1,\, i\ne j}^M\left(G_{ij}G'_{ji}+|\sigma_2|^2G_{ij}G'_{ij}\right)+\frac{\sigma_4-1}{M}\sum_{i=1}^M G_{ii}G'_{ii} \\
&=\frac{1}{\phi N}\sum_{i,j=1,\, i\ne j}^M\left(G_{ij}G'_{ji}+|\sigma_2|^2G_{ij}G'_{ij}\right)+(\sigma_4-1)m_{\phi^{-1}}(z)m_{\phi^{-1}}(z')\\
&\quad+\mathcal{O}_\prec\left(\frac{1}{|zz'|^\frac{1}{2}}\left(\frac{1}{\sqrt{N\eta}}+\frac{1}{\sqrt{N\eta'}}+\frac{1}{N\sqrt{\eta \eta'}}\right)\right),\numberthis \label{pair}
\end{align*} where we used that $\mathbb{E}{\bf x}_i^2=\mathbb{E}\widetilde{X}_{i1}^2=\frac{\sigma_2}{\sqrt{MN}}$ and $\mathbb{E}|{\bf x}_i|^4=\mathbb{E}|\widetilde{X}_{i1}|^4=\frac{\sigma_4}{MN}$ for each $i=1,\dots,M$. In the last step we also used \eqref{entrywise}.

To continue with the study of the fluctuation term we need to find an expression for $\frac{1}{\phi N}\sum_{i,j=1,\, i\ne j}^M G_{ij}G'_{ji}$ and $\frac{1}{\phi N}\sum_{i,j=1,\, i\ne j}^M G_{ij}G_{ij}'$ in terms of $m_\phi$ and $m_{\phi^{-1}}$.

\begin{lemma}
\label{selfconeq}
For $z=x+i\eta,\, z'=x'+i\eta'$, $\eta,\eta'>\eta_0$, with $|z-\sqrt{\phi}|\le 10$ and $|z'-\sqrt{\phi}|\le 10$, it holds \begin{equation}
\label{selftrace}
\frac{1}{\phi N}\sum_{i,j=1\atop i\ne j}^M G_{ij}G'_{ji}=\frac{zz'm_\phi(z)m_\phi(z')m_{\phi^{-1}}(z)^2m_{\phi^{-1}}(z')^2}{1-zz'm_\phi(z)m_\phi(z')m_{\phi^{-1}}(z)m_{\phi^{-1}}(z')}+\mathcal{O}_\prec\left(\frac{\bf{\Psi}}{|zz'|^\frac{1}{2}}\right),
\end{equation} 
\begin{equation}
\label{selfnontrace}
\frac{1}{\phi N}\sum_{i,j=1\atop i\ne j}^M G_{ij}G_{ij}'=\frac{|\sigma_2|^2zz'm_\phi(z)m_\phi(z')m_{\phi^{-1}}(z)^2m_{\phi^{-1}}(z')^2}{1-|\sigma_2|^2zz'm_\phi(z)m_\phi(z')m_{\phi^{-1}}(z)m_{\phi^{-1}}(z')}+\mathcal{O}_\prec\left(\frac{\bf{\Psi}}{|zz'|^\frac{1}{2}}\right),
\end{equation} where \[\bf{\Psi}:=\text{$\frac{1}{\eta+\eta'}\left(\frac{1}{\sqrt{N\eta\eta'^2}}+\frac{1}{\sqrt{N\eta^2\eta'}}+\frac{1}{N\eta\eta'}\right)$}.\]
\proof
\normalfont
To prove this lemma we change our point of view and we study the linearized problem. We remark that \eqref{selftrace}, being a tracial quantity, could still be analyzed without linearization, but \eqref{selfnontrace} cannot. For brevity we use the proof with linearization for both cases.

Let the $[(N-1)+M]\times [(N-1)+M]$ matrix 
$\mathcal{H}$ be defined as 
\begin{equation}
\label{big}
\mathcal{H}:=\left( \begin{matrix} 0 & X^* \\
 X & 0 \end{matrix}\right).
\end{equation} We introduced this bigger 
matrix $\mathcal{H}$ to study $W$, since $\mathcal{H}$ has the advantage that all nonzero elements are i.i.d. random variables (modulo symmetry) and it carries all information on the matrices $W=X^*X$ and $XX^*$ we are studying. Indeed, $\mathcal{H}^2$ with diagonal blocks $X^*X$ and $XX^*$ has the same non zero spectrum as $W$ (with double multiplicity).

To prove \eqref{selftrace} we define the resolvents \begin{equation}
\label{newresolvents}
\mathcal{G}(z):=(\mathcal{H}^2-z)^{-1}\,\,\, \text{and} \,\,\, \mathfrak{G}(\zeta):=(\mathcal{H}-\zeta)^{-1}.
\end{equation} Note that \begin{equation}
\label{reseq}
\mathcal{G}(z)=\frac{1}{2\sqrt{z}}\cdot\left(\frac{1}{\mathcal{H}-\sqrt{z}}-\frac{1}{\mathcal{H}+\sqrt{z}}\right)=\frac{1}{2\sqrt{z}}\cdot\left(\mathfrak{G}(\sqrt{z})-\mathfrak{G}(-\sqrt{z})\right),
\end{equation} where we chose the branch of $\sqrt{z}$ which lies in $\mathbb{H}$.

In the following we state some fundamental properties of the Gram matrix $\mathcal{H}$ and of its resolvent $\mathfrak{G}$ (for a detailed description see \cite{alt2017singularities} and \cite{alt2017local}). Let $m_1,m_2:\mathbb{H}\to \mathbb{H}$ be the unique solutions of the system \begin{equation}
\label{supers}
\begin{cases} -\frac{1}{m_1}=\zeta+\phi^\frac{1}{2}m_2, \\
-\frac{1}{m_2}=\zeta+\phi^{-\frac{1}{2}}m_1.
\end{cases}
\end{equation} Then, for each $\zeta\in\mathbb{H}$ (see \cite{alt2017local}) we have \begin{equation}
\label{entry}
\begin{split}
&|\mathfrak{G}_{ij}(\zeta)-\delta_{ij}m_1(\zeta)|\prec \frac{1}{\sqrt{N\Im\zeta}},\qquad i,j=2,\dots N, \\
&|\mathfrak{G}_{ij}(\zeta)-\delta_{ij}m_2(\zeta)|\prec \frac{1}{\sqrt{N\Im\zeta}}, \qquad i,j=N+1,\dots,N+M.
\end{split}
\end{equation} 

Notice that if $z=x+i\eta$ is such that $\zeta^2=z$ then $\frac{1}{\sqrt{N\Im \zeta}}\lesssim \frac{1}{\sqrt{N\eta}}$. Indeed, $\Im\zeta=\frac{\eta}{\Re \zeta}\gtrsim \eta$, since $|\zeta|\lesssim 1$ under the hypothesis $|z-\sqrt{\phi}|\le 10$ and \eqref{limitN}. Hence all along the proof we will estimate the error terms only in terms of $\eta$. We will use $\zeta$ as the argument of the resolvent $\mathfrak{G}$, with $\zeta=\sqrt{z}$.

In particular $m_1$ and $m_2$ are Stieltjes transforms of  symmetric probability measures on $\mathbb{R}$, whose support is contained in $[-2\phi^\frac{1}{4},2\phi^{-\frac{1}{4}}]$ (see Theorem 2.1 in \cite{ajanki2015quadratic}). Furthermore, we have that \begin{equation}
\label{tra}
m_\phi(z)=\frac{m_1(\zeta)}{\zeta}, \,\,\,m_{\phi^{-1}}(z)=\frac{m_2(\zeta)}{\zeta}
\end{equation} and they are related in the following way: \begin{equation}
\begin{split}
\begin{cases}
\label{win}
&-\frac{1}{m_\phi(z)}=z+z\phi^\frac{1}{2}m_{\phi^{-1}}(z)\\
&-\frac{1}{m_{\phi^{-1}}(z)}=z+z\phi^{-\frac{1}{2}}m_\phi(z).
\end{cases}
\end{split}
\end{equation} By \eqref{tra}, using that an analogue of \eqref{goodbound} holds substituting $\phi$ with $\phi^{-1}$ (see proof of Lemma \ref{sbound} in Appendix A), we have that \begin{equation}
\label{gm}
|\phi^{-\frac{1}{4}}m_1(z)|\le 1-c\eta,\,\,\,\,\, |\phi^\frac{1}{4}m_2(z)|\le 1-c\eta.
\end{equation}

Next, we use a resolvent expansion to  express the resolvents of $\mathcal{H}$ and $\mathcal{H}^2$ in terms of resolvents of their minors. For each $T\subset \{2,\dots, N+M\}$ we define \begin{equation}
\mathcal{G}^{[T]}(z):=\left(\left(\mathcal{H}^{[T]}\right)^2-z\right)^{-1}\,\,\, \text{and} \,\,\, \mathfrak{G}^{[T]}(\zeta):=(\mathcal{H}^{[T]}-\zeta)^{-1},
\end{equation} where $\mathcal{H}^{[T]}$ is the matrix $\mathcal{H}$ with the rows and columns labeled with $T$ set to zero:
\begin{equation}
\label{deleteh}
\left(\mathcal{H}^{[T]}\right)_{ij}:= \mathbb{1}(i\notin T)\mathbb{1} ( j\notin T) \mathcal{H}_{ij}.
\end{equation}
Let $\gamma_{ij}$ denote the entries of the matrix $\mathcal{H}$, i.e. $\gamma_{ij}=X_{ij}$ for $i=N+1,\dots N+M$, $j=2,\dots,N$, $\gamma_{ij}=\overline{\gamma_{ji}}$ for $i=2,\dots, N$, $j=N+1,\dots N+M$ and $\gamma_{ij}=0$ otherwise. From now on we abandon the convention in Remark \ref{ind} about Greek letters for columns indices and we use only $i,j,k,\dots$ We use the one sided expansion for the resolvent of $\mathcal{H}$, i.e. for each $i\ne j$ we have \begin{equation}
\label{oneside}
\mathfrak{G}_{ij}=-\mathfrak{G}_{ii}\sum_{k=2 \atop k\ne j}^{N+M}\mathfrak{G}_{ik}^{[j]}\mathcal{\gamma}_{kj}.
\end{equation} Notice that here $\mathfrak{G}_{ik}^{[j]}$ is independent of $\gamma_{kj}$ since $\mathcal{H}$ has independent elements. 

By the definition of $\mathcal{H}^2$ and \eqref{reseq}, using the identification $\zeta=\sqrt{z}$ choosing the branch of $\sqrt{z}$ which lies in $\mathbb{H}$, it follows that \begin{equation}
\label{relation}
\begin{split}\frac{1}{N}\sum_{i,j=1\atop i\ne j}^M G_{ij}G_{ji}'&=\frac{1}{N}\sum_{i,j=N+1\atop i\ne j}^{N+M}\mathcal{G}_{ij}(z)\mathcal{G}_{ji}(z')\\
&=\frac{1}{N}\sum_{i,j=N+1\atop i\ne j}^{N+M}\frac{1}{4\zeta\zeta'}\left(\mathfrak{G}(\zeta)_{ij}\mathfrak{G}(\zeta')_{ji}-\mathfrak{G}(\zeta)_{ij}\mathfrak{G}(-\zeta')_{ji}\right)\\
&\quad+\frac{1}{N}\sum_{i,j=N+1\atop i\ne j}^{N+M}\frac{1}{4\zeta\zeta'}\left(\mathfrak{G}(-\zeta)_{ij}\mathfrak{G}(-\zeta')_{ji}-\mathfrak{G}(-\zeta)_{ij}\mathfrak{G}(\zeta')_{ji}\right).\end{split}
\end{equation} 

We introduce the shorthand notation $\mathfrak{G}_{ij}:=\mathfrak{G}_{ij}(\zeta)$, $\mathfrak{G}_{ij}':=\mathfrak{G}_{ij}(\zeta')$. By \eqref{entry}, for any $i,j,k$ all distinct, it holds \begin{equation}
\label{ri}
\mathfrak{G}_{ik}=\mathfrak{G}_{ik}^{[j]}+\frac{\mathfrak{G}_{ij}\mathfrak{G}_{jk}}{\mathfrak{G}_{jj}}=\mathfrak{G}_{ik}^{[j]}+\mathcal{O}_\prec\left(\frac{1}{N\eta}\right).
\end{equation}
We now derive a self consistent equation for $\sum_{i\ne j}\mathfrak{G}_{ij}\mathfrak{G}_{ji}'$, that is the first term in the second equality of \eqref{relation}.

For this purpose, we start proving that $\sum_{i\ne j}\mathfrak{G}_{ij}\mathfrak{G}_{ji}'$ is close to  $\sum_{i\ne j} \mathbb{E}_j \mathfrak{G}_{ij}\mathfrak{G}_{ji}'$
where $\mathbb{E}_j(\cdot):=\mathbb{E}(\cdot|\mathcal{H}^{[j]})$ denotes the conditional expectation with respect to the matrix $\mathcal{H}^{[j]}$. 
This result is a  special case of the fluctuation  averaging analysis presented in \cite{erdHos2013averaging}, 
in fact its very elementary version given in Proposition 6.1 of \cite{erdHos2013averaging} suffices.  
No other input from the technically involved paper \cite{erdHos2013averaging} is used for the proof of \eqref{me}.
More precisely, for any fixed $i$, we have the bound
\begin{equation}
\label{me}
\frac{1}{N}\sum_{j=N+1\atop j\ne i}^{N+M}(1-\mathbb{E}_j)\mathfrak{G}_{ij}\mathfrak{G}_{ji}'=\mathcal{O}_\prec\left(\frac{1}{\sqrt{N\eta}}\frac{1}{\sqrt{N\eta'}}\left(\frac{1}{\sqrt{N\eta}}+\frac{1}{\sqrt{N\eta'}}\right)\right).
\end{equation}
In particular, \eqref{me} shows that the operator $(1-\mathbb{E}_j)$ reduces the naive size of $\frac{1}{N} \sum_{i\ne j}\mathfrak{G}_{ij}\mathfrak{G}_{ji}'$
coming from \eqref{entry} by an additional factor $1/\sqrt{N\eta}+1/\sqrt{N\eta'}$. Indeed, by \cite[Eq. (4.5)]{erdHos2013averaging}, the left hand side of \eqref{me} is exactly the left hand side of \cite[Eq. 6.1]{erdHos2013averaging} after the associations ${\bf \emph{a}}=(i)$, ${\bf \mu}=(j)$, $w({\bf \emph{a}})=w(i)=N^{-1}$, $F=\{j\}$ and $\Delta$ being the graph of degree deg$(\Delta)=2$ corresponding to $\mathfrak{G}_{ij}\mathfrak{G}_{ji}'$. Now we explain the single modification in the  proof of Proposition 6.1 in \cite{erdHos2013averaging} that leads to \eqref{me}.

We recall that the main strategy in the proof of Proposition 6.1 in \cite{erdHos2013averaging} is to compute the $p$-th moment 
of the sum $\sum_j (1-\mathbb{E}_j)\mathfrak{G}_{ij}\mathfrak{G}_{ji}'$. Expanding the $p$-th power yields a $p$-fold summation $\sum_{j_1, j_2, \ldots, j_p}$. For any fixed choice of these indices, we 
 successively  expand the resolvent entries  as much as possible, in order to create factors  partially independent of each other
  using the resolvent expansion \eqref{ri} for terms of the form $\mathfrak{G}_{ik}^{[T]}$, with $i,k\notin T$, and its analogues for $1/\mathfrak{G}_{ii}^{[T]}$ from \cite[Eq. (3.13)]{erdHos2013averaging}.
 Here the set $T$ is a subset of the actual summation indices $j_1, j_2, \ldots , j_p$. After taking the expectation and using that $\mathbb{E}(1-\mathbb{E}_j)=0$, 
 a simple power counting shows that only those terms remain nonzero that have many resolvent factors.
 Then, after that each factor is expanded as described above, we use the bound $|\mathfrak{G}_{ij}(z)|\le 1/\sqrt{N\Im z}$, given by the local law in \eqref{entry} for $i\ne j$. In particular, in the proof of Proposition 6.1 in \cite{erdHos2013averaging} the resolvent expansions and the bounds given by the local law are used only for single resolvent entries. Hence, 
 the proof of Proposition 6.1 \cite{erdHos2013averaging}  works verbatim for our case when different spectral parameter are considered, just in the estimates the different $\eta$'s have to be 
 carried.
  As a consequence, the error term in the r.h.s of \eqref{me}, in contrast to its analogue in \cite[Eq. (6.1)]{erdHos2013averaging}, contains both $\eta$ and $\eta'$, i.e. the error term is of the form $1/\sqrt{N^3\eta^2\eta'}+1\sqrt{N^3\eta\eta'^2}$.
 
By \eqref{me}, \eqref{oneside} and the local laws in \eqref{entry} we get \begin{align*}
\frac{1}{N}\sum_{i,j=N+1\atop j\ne i}^{N+M}\mathfrak{G}_{ij}\mathfrak{G}_{ji}'&=\frac{1}{N}m_2(\zeta)m_2(\zeta')\sum_{i,j=N+1\atop j\ne i}^{N+M}\mathbb{E}_j\left(\sum_{k=2\atop k\ne j}^{N+M}\mathfrak{G}_{ik}^{[j]}\gamma_{kj}\right)\left(\sum_{l=2\atop l\ne j}^{N+M}\gamma_{jl}\mathfrak{G}_{li}'^{[j]}\right)\\
&\quad +\mathcal{O}_\prec\left((\eta+\eta')\bf{\Psi}\right)\numberthis\label{supereqres1}\\
&=\frac{1}{N\sqrt{MN}}m_2(\zeta)m_2(\zeta')\sum_{i,j=N+1 \atop j\ne i}^{N+M}\sum_{k=2}^N\mathfrak{G}_{ik}^{[j]}\mathfrak{G}_{ki}'^{[j]}+\mathcal{O}_\prec\left((\eta+\eta')\bf{\Psi}\right).
\end{align*}
Note that we used \eqref{tra} and \eqref{gm} to estimate the error terms. Using \eqref{ri} the resolvent expansion in \eqref{oneside} and fluctuation averaging \eqref{me} again, \eqref{supereqres1} becomes
\begin{equation}
\label{supereqres}
\begin{split} \frac{1}{N}\sum_{i,j=N+1\atop j\ne i}^{N+M}\mathfrak{G}_{ij}\mathfrak{G}_{ji}'&=\frac{\sqrt{\phi}}{N}m_2(\zeta)m_2(\zeta')\sum_{i=N+1}^{N+M}\sum_{k=2}^N \mathfrak{G}_{ik}\mathfrak{G}_{ki}' +\mathcal{O}_\prec\left((\eta+\eta')\bf{\Psi}\right)\\
&=\frac{\sqrt{\phi}}{N}m_2(\zeta)m_2(\zeta')\sum_{i=N+1}^{N+M}\sum_{k=2}^N \mathbb{E}_k\mathfrak{G}_{ik}\mathfrak{G}_{ki}' +\mathcal{O}_\prec\left((\eta+\eta')\bf{\Psi}\right)\\
&=\frac{1}{N}m_1(\zeta)m_1(\zeta')m_2(\zeta)m_2(\zeta')\sum_{i,p=N+1\atop p\ne i}^{N+M}\mathfrak{G}_{ip}\mathfrak{G}_{pi}'\\
&\quad+\phi m_1(\zeta)m_1(\zeta')m_2(\zeta)^2m_2(\zeta')^2+\mathcal{O}_\prec\left((\eta+\eta')\bf{\Psi}\right).
\end{split}
\end{equation}

Solving this equation, we conclude that \begin{equation}
\label{primo}
\frac{1}{N}\sum_{i,j=N+1\atop i\ne j}^{N+M} \mathfrak{G}_{ij}\mathfrak{G}_{ji}'=\frac{\phi m_1(\zeta)m_1(\zeta')m_2(\zeta)^2m_2(\zeta')^2}{1-m_1(\zeta)m_1(\zeta')m_2(\zeta)m_2(\zeta')}+\mathcal{O}_\prec\left(\bf{\Psi}\right).
\end{equation} In estimating the error term we used a lower bound for the denominator. Indeed, using \eqref{tra} and \eqref{gm}, we have that \begin{equation}
\label{ser}
|1-m_1(\zeta)m_1(\zeta')m_2(\zeta)m_2(\zeta')|\ge 1-|m_1(\zeta)m_1(\zeta')m_2(\zeta)m_2(\zeta')|\gtrsim(\eta+\eta').
\end{equation}

Notice that in the right hand side of \eqref{primo} the deterministic term depends only on $m_1$ and $m_2$. Moreover, using the notation $\widehat{\mathfrak{G}}(\zeta):=(-\mathcal{H}-\zeta)^{-1}$ and that $m_1$ and $m_2$ are Stieltjes transforms of symmetric distributions, by \eqref{entry} we have that \begin{align}
&\left|\widehat{\mathfrak{G}}_{ij}(\zeta)-\delta_{ij}m_1(\zeta)\right|\prec \frac{1}{\sqrt{N\eta}},\,\,\,\,\,i,j=2,\dots, ,N\label{secondo}\\
&\left|\widehat{\mathfrak{G}}_{ij}(\zeta)-\delta_{ij}m_2(\zeta)\right|\prec \frac{1}{\sqrt{N\eta}},\,\,\,\,\, i,j=N+1,\dots,N+M\label{terzo}.
\end{align} In \eqref{secondo} and \eqref{terzo} we used that $\Im\zeta\gtrsim\eta$. This means that the leading order deterministic term of each term in \eqref{relation} is exactly the same. Hence, combining \eqref{relation}, \eqref{primo} and \eqref{tra} we conclude \eqref{selftrace}. The proof of \eqref{selfnontrace} is analogous.
\endproof
\end{lemma} 

Before proceeding, we recall that $f_\mathbb{C}(z)$ is supported in $[a,b]\times [-10,10]$, where $a=\gamma_--\delta_*$, $b=\gamma_++\delta_*$  and $\gamma_-$, $\gamma_+$ are the spectral edges. Furthermore, we recall that, by \eqref{wwin}, $w_\phi=-zm_\phi(z)m_{\phi^{-1}}(z)$, where $w_\phi(z)$ is the Stieltjes transform of the Wigner semicircle law centered at $\phi^\frac{1}{2}+\phi^{-\frac{1}{2}}$, hence $w_\phi(z)$ is a solution of the self consistent equation \eqref{stw}.

We now plug \eqref{pair}--\eqref{selfnontrace} into the integral in \eqref{var56}. Integrating the error terms in \eqref{pair}--\eqref{selfnontrace} and using that $|g(z)|\le C|z|^\frac{1}{2}$ (see \eqref{goodbound} and \eqref{fsmall}) we get an error term of the magnitude $N^{-\frac{7}{6}}$. The denominators in \eqref{selftrace} and \eqref{selfnontrace} are expanded into geometric series whose convergence follows from \eqref{tra} and \eqref{ser}. Hence, using \eqref{wwin}, we conclude that if $\sigma_2=0$ then \eqref{var56} assumes the following form\begin{align*}
F_N'&=\frac{1}{2N}\Re\iint_a^b\iint_{\eta_0}^{10}\Big[g(z)g(\overline{z}')m_{\phi^{-1}}(z)m_{\phi^{-1}}(\overline{z}')\sum_{k\ge 1}\big[w_\phi(z)w_\phi(\overline{z}')\big]^k\\
&\quad-g(z)g(z')m_{\phi^{-1}}(z)m_{\phi^{-1}}(z')\sum_{k\ge 1}\big[w_\phi(z)w_\phi(z')\big]^k\Big]d\eta d\eta' dxdx'\\
&\quad+\frac{\sigma_4-1}{N}\left(\Im \int_a^b\int_{\eta_0}^{10} g(z)m_{\phi^{-1}}(z)\,d\eta dx\right)^2+\mathcal{O}_\prec\left(N^{-\frac{7}{6}}\right)\\
&=\frac{1}{N}\sum_{k\ge 1}\left(\Im\int_a^b\int_{\eta_0}^{10} g(z)m_{\phi^{-1}}(z)w_\phi(z)^k\,d\eta dx\right)^2\numberthis\label{A}\\
&\quad+\frac{\sigma_4-1}{N}\left(\Im \int_a^b\int_{\eta_0}^{10} g(z)m_{\phi^{-1}}(z)\,d\eta dx\right)^2+\mathcal{O}_\prec\left(N^{-\frac{7}{6}}\right).\\
\end{align*} Substituting the expression of $g$ (see \eqref{short}) in \eqref{A} we have \begin{equation}
\label{Z}
\begin{split}
F_N'&=\frac{1}{N}\sum_{k\ge 2}\left(\frac{2}{\pi}\Im\int_a^b\int_{\eta_0}^{10} w_\phi(z)^k\partial_\eta\partial_{\overline{z}}f_\mathbb{C}(z)\,d\eta dx\right)^2\\
&\quad+\frac{\sigma_4-1}{N}\left(\frac{2}{\pi}\Im \int_a^b\int_{\eta_0}^{10} w_\phi(z)\partial_\eta\partial_{\overline{z}}f_\mathbb{C}(z)\,d\eta dx\right)^2+\mathcal{O}_\prec\left(N^{-\frac{7}{6}}\right).
\end{split}\end{equation} 
We start computing the last integral in \eqref{Z}: 
\[
\left(\frac{2}{\pi}\Im \int_a^b\int_{\eta_0}^{10} w_\phi(z)\partial_\eta\partial_{\overline{z}}f_\mathbb{C}(z)\,d\eta dx\right)^2
=\left(\frac{1}{\pi}\Im \int_a^b w_\phi(x)f'(x)\, dx\right)^2+\mathcal{O}_\prec\left(N^{-\frac{1}{6}}\right),
\]

where we used Lemma \ref{sto} and \begin{equation}
\label{fpsmall}
\frac{\partial_\eta f_\mathbb{C}(z_0)}{i}=\partial_x f_\mathbb{C}(z_0)+\mathcal{O}(\eta_0)=f'(x)+\mathcal{O}_\prec(\eta_0),
\end{equation} where $z_0=x+i\eta_0$. Furthermore, using Lemma \ref{sto} and \eqref{fpsmall} once more, we have \begin{align*}
&\frac{1}{N}\sum_{k\ge 2}\left(\frac{2}{\pi}\Im\int_a^b\int_{\eta_0}^{10} w_\phi(z)^k\partial_\eta\partial_{\overline{z}} f_\mathbb{C}(z)\,d\eta  dx\right)^2\numberthis \label{C}\\
&=\frac{1}{N}\sum_{k\ge 2}\left(\frac{1}{\pi}\Im\int_a^b w_\phi(z_0)^k\frac{\partial_\eta f_\mathbb{C}(z_0)}{i}\, dx\right)^2+\mathcal{O}_\prec\left(N^{-\frac{7}{6}}\right)\\
&=\frac{1}{N}\sum_{k\ge 0}\left(\frac{1}{\pi}\Im\int_a^bw_\phi(z_0)^k\frac{\partial_\eta f_\mathbb{C}(z_0)}{i}\,dx\right)^2 -\frac{1}{N}\left(\frac{1}{\pi}\Im\int_a^bw_\phi(z_0)\frac{\partial_\eta f_\mathbb{C}(z_0)}{i}\, dx\right)^2 \\
&\quad-\frac{1}{N}\left(\frac{1}{\pi}\Im\int_a^b\frac{\partial_\eta f_\mathbb{C}(z_0)}{i}\, dx\right)^2+\mathcal{O}_\prec\left(N^{-\frac{7}{6}}\right)\\
&=\frac{1}{N}\sum_{k\ge 0}\left(\frac{1}{\pi}\Im\int_a^bw_\phi(z_0)^k\frac{\partial_\eta f_\mathbb{C}(z_0)}{i}\, dx\right)^2-\left(\frac{1}{\pi}\Im\int_a^b w_\phi(x)f'(x)\, dx\right)^2+\mathcal{O}_\prec\left(N^{-\frac{7}{6}}\right).
\end{align*} In the last equality we used that $\Im \frac{\partial_\eta f_\mathbb{C}(z_0)}{i}=\mathcal{O}_\prec(\eta_0)$ by \eqref{fpsmall}. We want to use the same approximation in the first integral as well. However, the geometric series converges only slowly, so we need to ensure summability. The following lemma prepares us for that (see Lemma 3.7 in \cite{erdHos2017fluctuations}).

\begin{lemma}
\label{usbound}
There exists an $N$-independent constant $C>0$ such that for $z_0=x+i\eta_0$ and $z_0'=x'+i\eta_0$, with $0<\eta_0\le \frac{1}{2}$, it holds \begin{equation}
\label{bbb}
\iint_a^b\frac{dxdx'}{|1-w_\phi (z_0)w_\phi(\overline{z_0}')|}+\iint_a^b\frac{dxdx'}{|1-w_\phi (z_0)w_\phi(z_0')|}\le C|\log\eta_0|.
\end{equation}
\end{lemma}

Combining \eqref{Z}-\eqref{C} and Lemma \ref{usbound}, using \eqref{complex} again, we conclude that \begin{align*}
F_N'&=\frac{1}{2N\pi^2}\Re \iint_a^b\left(\frac{1}{1-w_\phi (z_0)w_\phi(\overline{z_0}')}-\frac{1}{1-w_\phi (z_0)w_\phi(z_0')}\right)f'(x)f'(x')\, dxdx'\\
&\quad+\frac{\sigma_4-2}{N}\left(\frac{1}{\pi}\Im \int_a^b w_\phi(x)f'(x)\, dx\right)^2+\mathcal{O}\left(N^{-\frac{7}{6}}\right). \numberthis\label{E}
\end{align*} 

After some computations  using \eqref{stw} we have that \begin{equation}
\label{F}
\begin{split}
\Re &\left(\frac{1}{1-w_\phi (z_0)w_\phi(\overline{z_0}')}-\frac{1}{1-w_\phi (z_0)w_\phi(z_0')}\right)\\
&=\Re \left( \frac{2i\Im w_\phi(z_0')}{\phi^\frac{1}{2}+\phi^{-\frac{1}{2}}-z_0-2\Re w_\phi(z_0')-w_\phi(z_0')(|w_\phi(z_0')|^2-1)} \right).
\end{split}\end{equation} For small $\eta_0$ and $(x,x')$ outside the square $[\gamma_-,\gamma_+]^2$ the integral of \eqref{F} is negligible. Indeed, outside $[\gamma_- ,\gamma_+]^2$ we have that $1-|w_\phi(z)|^2\asymp \sqrt{\kappa_x+\eta}$ by Lemma~3.6 in \cite{bloemendal2016principal}, where $\kappa_x=\min\{|\gamma_+-x|,|\gamma_--x|\}$.

For $(x,x')\in [\gamma_-,\gamma_+]^2$ and small $\eta_0$ we have \begin{equation}
\label{W}
\begin{split}
\Re &\left( \frac{2i\Im w_\phi(z_0')}{\phi^\frac{1}{2}+\phi^{-\frac{1}{2}}-z_0-2\Re w_\phi(z_0')-w_\phi(z_0')(|w_\phi(z_0')|^2-1)} \right)\\
&\quad =\frac{\eta_0\sqrt{(x'-\gamma_-)(\gamma_+-x')}}{(x-x')^2+\eta_0^2}+\mathcal{O}_\prec(\eta_0).
\end{split}\end{equation} The expression $\frac{\eta_0}{(x-x')^2+\eta_0^2}$ acts like $\pi\delta(x'-x)$ for small $\eta_0$, hence for each $h\in L^2$ \[\lim_{\eta\to 0}\int_\mathbb{R}\frac{\eta}{(x-x')^2+\eta^2}h(x')\,dx'=\pi h(x)\] in $L^2$-sense. Working out an effective error term for $h\in H^1$ and using the explicit expression in \eqref{W}, by \eqref{E}, we conclude that \[ \begin{split}F_N'&=\frac{1}{2\pi N}\int_{\gamma_-}^{\gamma^+} f'(x)^2\sqrt{(x-\gamma_-)(\gamma_+-x)}\,dx\\
&\quad+\frac{\sigma_4-2}{N}\left(\frac{1}{\pi}\int_{\gamma_-}^{\gamma_+} \frac{1}{2}f'(x)\sqrt{(x-\gamma_-)(\gamma_+-x)}\, dx\right)^2+\mathcal{O}\left(N^{-\frac{7}{6}}\right).\end{split}\] This computation gives the explicit expression of $V_f$ in \eqref{var2} for $\sigma_2=0$.

When $\sigma_2\ne 0$ we have to consider \eqref{selfnontrace} and so, using a similar analysis, we have to add the following term in the expression of $F_N'$ in \eqref{A}\begin{equation}
\label{arar}
\frac{1}{2N\pi^2}\Re \iint_\mathbb{R}f'(x)f'(x')\left(\frac{|\sigma_2|^2w_\phi(z_0)^2 w_\phi(\overline{z_0}')^2}{1-|\sigma_2|^2w_\phi(z_0)w_\phi(\overline{z_0}')}-\frac{|\sigma_2|^2w_\phi(z_0)^2w_\phi(z_0')^2}{1-|\sigma_2|^2w_\phi(z_0)w_\phi(z_0')}\right)\,dxdx'.
\end{equation}

For the special case $|\sigma_2|=1$ the expressions in \eqref{selftrace} and \eqref{selfnontrace} are exactly the same, hence we define $V_{\sigma_2}:=V_{f,1}$. This holds true in particular for the case $X\in \mathbb{R}^{M\times (N-1)}$ when $\sigma
_2=1$ automatically.

If $|\sigma_2|<1$, instead, we define $V_{\sigma_2}$ in the following way 
\begin{equation}
\label{G}
V_{\sigma_2}:=\frac{1}{2\pi^2}\Re \iint_\mathbb{R}f'(x)f'(x')\left(\frac{|\sigma_2|^2w_\phi(x)^2 \overline{w_\phi(x')}^2}{1-|\sigma_2|^2w_\phi(x)
\overline{w_\phi(x')}}-\frac{|\sigma_2|^2w_\phi(x)^2w_\phi(x')^2}{1-|\sigma_2|^2w_\phi(x)w_\phi(x')}\right)dxdx',
\end{equation}
 that is close to \eqref{arar} by an $\mathcal{O}(\eta_0)$ error using that $|w_\phi(z_0)-w_\phi(x)|\lesssim \eta_0[(x-\gamma_-)(\gamma_+-x)]^{-\frac{1}{2}}$ and $|1-|\sigma_2|^2w_\phi(x)\overline{w_\phi(x')}|\ge 1-|\sigma_2|^2$. Notice that from \eqref{G} easily follows that $V_{\sigma_2}\ge 0$. Indeed \[V_{\sigma_2}=\sum_{k\ge0}\left(\frac{1}{\pi}\Im\int_a^bf'(x)(|\sigma_2|w_\phi(x))^{k+2}\,dx\right)^2.\]

\section{Computation of the higher order moments of $F_N$}
In this section we compute the higher order moments of \[F_N=-\frac{1}{\sqrt{N}}\Im\int_\mathbb{R}\int_{\eta_0}^{10} g(z)A(z)\,d\eta dx+\mathcal{O}_\prec(\eta_0),\] where $g(z)$ and $A(z)$ are defined in \eqref{short}. We remark that for the proof of the normality of $F_N$ it would be sufficient to show that the quadratic form $\langle {\bf x}, G(z) {\bf x}\rangle$ has a 
Gaussian fluctuation conditioned on $G$ and then separately show that the quadratic variation of $G$ is negligible. Here we follow a more robust path that gives an
effective control on all higher moments as well without essentially no extra effort since the fluctuation averaging mechanism used already in the proof of Lemma~\ref{selfconeq}
directly  extends to higher moments. Thus, using a similar approach to the one we used to compute the variance of $F_N$, we start computing\[\mathbb{E}[A(z_1)\dots A(z_k)]\] for any $k\in\mathbb{N}$ and $z_l\in\mathbb{C}\setminus\mathbb{R}$, with $l=1,\dots, k$. We recall that $\mathbb{E}_1:= \mathbb{E}(\cdot|X)$ is the conditional expectation conditioned on the matrix $X$. This leads to products of cyclic expressions of the form $G_{j_1j_2}G_{j_2j_3}\dots G_{j_{k-1}j_k}$.

{\bf Notation}. A multiple summation with a star  $\sum_{j_1,\dots,j_k}^*$ indicates that the sum is performed over distinct indices.

In the following we prove that the leading order term of the $k$-th moment of $F_N$ is given by cycles of length two, hence cyclic products with at least three terms are actually of lower order:

\begin{lemma}
\label{cyc}
For closed cycles of length $k>2$ we have that \begin{equation}
\label{eq1}
N^{-\frac{k}{2}}\sum_{j_1,\dots, j_k=1}^{M\,\,*}\mathbb{E}_{j_1+N}\left(G_{j_1j_2}^{(1)}\dots G_{j_{k-1}j_k}^{(k-1)}G_{j_kj_1}^{(k)}\right)\prec \frac{|z_1\dots z_k|^{-\frac{1}{2}}}{(\max_a\eta_a)\sqrt{N\eta_1\dots\eta_k}}\sum_{a=1}^k\frac{1}{\sqrt{\eta_a}},
\end{equation} and for open cycles of any length $k>1$ we have that \begin{equation}
\label{eq2}
N^{-\frac{k+1}{2}}\sum_{j_1,\dots, j_k=1}^{M\,\, *} \mathbb{E}_{j_1+N}\left(G_{j_1j_2}^{(1)}\dots G_{j_{k-1}j_k}^{(k-1)}\right)\prec \frac{|z_1\dots z_k|^{-\frac{1}{2}}}{\sqrt{N \eta_1\dots\eta_{k-1}}}\sum_{a=1}^k\frac{1}{\sqrt{\eta_a}},           
\end{equation} where $G^{(l)}:=G(z_l)$, $z_l\in\mathbb{C}\setminus\mathbb{R}$ with $\eta_l=|\Im z_l|$ for $l=1,\dots,k$ and $\mathbb{E}_{j_1+N}:=\mathbb{E}(\cdot|\mathcal{H}^{[j_1+N]})$, with $\mathcal{H}^{[j_1+N]}$ defined in \eqref{deleteh}. Moreover, the same bounds hold true when any of the $G^{(l)}$ are replaced by their transposes or Hermitian conjugates.
\proof
\normalfont
The proof is similar to the proof of Lemma 4.1 in \cite{erdHos2017fluctuations}, so we will skip some details. However, an additional step in needed, see \eqref{eq7} later.

We start proving \eqref{eq1} for the case $X\in\mathbb{R}^{M\times (N-1)}$. We will actually prove that \[N^{-\frac{k}{2}}\sum_{j_1,\dots, j_k=1}^{M\,\, *} \mathbb{E}_{j_1+N}(G_{j_1j_2}^{(1)}\dots G_{j_{k-1}j_k}^{(k-1)}G_{j_kj_1}^{(k)})\lesssim \frac{N^\epsilon |z_1\dots z_k|^{-\frac{1}{2}}}{(\eta_1+\eta_k)\sqrt{N\eta_1,\dots,\eta_k}}\sum_{a=1}^k\frac{1}{\sqrt{\eta_a}},\] for any $\epsilon>0$, which implies \eqref{eq1} by the definition of $\prec$ in Definition \ref{sd}.

We use linearization again to express the resolvents $G^{(1)},\dots,G^{(k)}$ of the matrix $XX^*$ in terms of the resolvents $\mathfrak{G}^{(1)},\dots,\mathfrak{G}^{(k)}$ of the linearized matrix $\mathcal{H}$. 

\begin{equation}
\label{eq3}
\sum_{j_1,\dots, j_k=1}^{M\,\,*}\mathbb{E}_{j_1+N}\left(G_{j_1j_2}^{(1)}\dots G_{j_{k-1}j_k}^{(k-1)}G_{j_kj_1}^{(k)}\right)=\sum_{i_1,\dots, i_k=N+1}^{N+M\, *} \mathbb{E}_{i_1}\left(\mathcal{G}_{i_1i_2}^{(1)}\dots \mathcal{G}_{i_{k-1}i_k}^{(k-1)}\mathcal{G}_{i_ki_1}^{(k)}\right),
\end{equation} where $\mathcal{G}^{(l)}=(\mathcal{H}^2-z_l)^{-1}$, $i_m=j_m+N$. We write each $\mathcal{G}^{(l)}$ in the r.h.s. of \eqref{eq3} as \begin{equation}
\label{aaa}
\mathcal{G}(z_l)=\frac{1}{2\zeta_l}\cdot\left(\mathfrak{G}(\zeta_l)-\mathfrak{G}(-\zeta_l)\right),
\end{equation} with $\zeta_l^2=z_l$ (see \eqref{reseq}). We have to find a self consistent equation for each term in the right-hand side of \eqref{eq3} after rewriting it using \eqref{aaa}. We start with \[N^{-\frac{k}{2}}\sum_{i_1,\dots, i_k=N+1}^{N+M\, *} \mathbb{E}_{i_1}\left(\mathfrak{G}_{i_1i_2}^{(1)}\dots \mathfrak{G}_{i_{k-1}i_k}^{(k-1)}\mathfrak{G}_{i_ki_1}^{(k)}\right).\] Using the resolvent identity $\mathfrak{G}^{(1)}=\frac{1}{\zeta_1}[\mathcal{H}^{(1)}\mathfrak{G}^{(1)}-1]$ we get \begin{equation}
\label{eq4}
\begin{split}
N^{-\frac{k}{2}}&\sum_{i_1,\dots, i_k=N+1}^{N+M\, *} \mathbb{E}_{i_1}\left(\mathfrak{G}_{i_1i_2}^{(1)}\dots \mathfrak{G}_{i_{k-1}i_k}^{(k-1)}\mathfrak{G}_{i_ki_1}^{(k)}\right)\\
&=\frac{1}{N^\frac{k}{2}\zeta_1}\sum_{i_1,\dots, i_k=N+1}^{N+M\, *}\sum_{n=2}^{N+M} \mathbb{E}_{i_1}\left(\gamma_{i_1n}\mathfrak{G}_{ni_2}^{(1)}\dots \mathfrak{G}_{i_{k-1}i_k}^{(k-1)}\mathfrak{G}_{i_ki_1}^{(k)}\right),
\end{split}\end{equation} where $\gamma_{ij}$, with $i,j\in\{2,\dots,N+M\}$, are the entries of the big matrix $\mathcal{H}$.

We use the standard cumulant expansion \begin{equation}
\label{cumulant}
\begin{split}
\mathbb{E}hf(h)&=\mathbb{E}h\mathbb{E}f(h)+\mathbb{E}h^2\mathbb{E}f'(h)+\mathcal{O}\left(\mathbb{E}\left|h^3 \mathbb{1}(|h|>N^{\tau-\frac{1}{2}})\right|\Vert f''\rVert_\infty\right)\\
&\quad +\mathcal{O}\Big(\mathbb{E}|h|^3\sup_{|x|\le N^{\tau-\frac{1}{2}}}|f''(x)|\Big),
\end{split}\end{equation} where $f$ is any smooth function of a real random variable $h$, such that the expectations exist and $\tau>0$ is arbitrary (see \cite{khorunzhy1996asymptotic}). This yields \begin{align*}
 \mathbb{E}_{i_1}\left(\gamma_{i_1n}\mathfrak{G}_{ni_2}^{(1)}\dots \mathfrak{G}_{i_{k-1}i_k}^{(k-1)}\mathfrak{G}_{i_ki_1}^{(k)}\right)&=\frac{1}{\sqrt{MN}}\mathbb{E}_{i_1}\left(\frac{\partial\mathfrak{G}_{ni_2}^{(1)}}{\partial\gamma_{i_1n}}\mathfrak{G}_{i_2i_3}^{(2)}\dots\mathfrak{G}_{i_ki_1}^{(k)}\right)\numberthis\label{eq5}\\
&\quad+\frac{1}{\sqrt{MN}}\sum_{a=2}^k\mathbb{E}_{i_1}\left(\frac{\partial\mathfrak{G}_{i_ai_{a+1}}^{(a)}}{\partial \gamma_{i_1n}}\mathfrak{G}_{ni_2}^{(1)}\prod_{a\ne b=2}^k\mathfrak{G}_{i_bi_b+1}^{(b)}\right)+R,
\end{align*} where $i_{k+1}=i_1$ and $R$ is the error term resulting from the cumulant expansion.

Using the expression for the derivative of the resolvent \[\frac{\partial\mathfrak{G}_{ij}}{\partial\gamma_{kl}}=-\frac{\mathfrak{G}_{ik}\mathfrak{G}_{lj}+\mathfrak{G}_{il}\mathfrak{G}_{kj}}{1+\delta_{kl}}\] and the local law by \eqref{entry} for the resolvent of the Gram matrix $\mathcal{H}$, summing over $n$, the first term of the right hand side of \eqref{eq5} becomes \begin{equation}
\label{eq6}
\begin{split}
-\frac{1}{\sqrt{MN}} &\sum_{n=2}^N\left(\mathfrak{G}_{ni_1}^{(1)}\mathfrak{G}_{ni_2}^{(1)}+\mathfrak{G}_{nn}^{(1)}\mathfrak{G}_{i_1i_2}^{(1)}\right)\mathfrak{G}_{i_2i_3}^{(2)}\dots \mathfrak{G}_{i_ki_1}^{(k)}\\
&\quad =-\phi^{-\frac{1}{2}}m_1(\zeta_1)\mathfrak{G}_{i_1i_2}^{(1)}\dots \mathfrak{G}_{i_ki_1}^{(k)}+\mathcal{O}_\prec \left(\frac{1}{N^{\frac{k}{2}+\frac{1}{2}}\sqrt{\eta\eta_1}}\right),
\end{split}\end{equation} with $n\ne i_1, i_2$ and $\eta:=\eta_1\dots \eta_k$. If $n$ is equal to $i_1$ or $i_2$ we use the trivial bound. 
\\
Using the same computations of Lemma 4.1 in \cite{erdHos2017fluctuations}, if $a\ne k$ the second term of the right-hand side of \eqref{eq5} can be estimated by \[-\left(\mathfrak{G}_{i_ai_1}^{(a)}\mathfrak{G}_{ni_{a+1}}^{(a)}+\mathfrak{G}_{i_a n}^{(a)}\mathfrak{G}_{i_1i_{a+1}}^{(a)}\right)\mathfrak{G}_{n i_2}^{(1)}\prod_{a\ne b=2}^k \mathfrak{G}_{i_bi_{b+1}}^{(b)}\prec \frac{1}{N^{\frac{k}{2}}\sqrt{\eta\eta_a}}\] and if $n\notin\{i_1,\dots,i_k\}$ this bound can be improved to \[-\left(\mathfrak{G}_{i_ai_1}^{(a)}\mathfrak{G}_{ni_{a+1}}^{(a)}+\mathfrak{G}_{i_a n}^{(a)}\mathfrak{G}_{i_1i_{a+1}}^{(a)}\right)\mathfrak{G}_{n i_2}^{(1)}\prod_{a\ne b=2}^k \mathfrak{G}_{i_bi_{b+1}}^{(b)}\prec \frac{1}{N^{\frac{k}{2}+\frac{1}{2}}\sqrt{\eta\eta_a}}.\]

 Finally, for the case $a=k$ we have \[ -\left(\mathfrak{G}_{i_ki_1}^{(k)}\mathfrak{G}_{ni_1}^{(k)}+\mathfrak{G}_{i_k n}^{(k)}\mathfrak{G}_{i_1i_1}^{(k)}\right)\mathfrak{G}_{n i_2}^{(1)}\dots \mathfrak{G}_{i_{k-1}i_k}^{(k-1)}.\] Here an additional argument is needed compared to \cite{erdHos2017fluctuations}. To get a similar expression to \eqref{eq6} we need to have that all the indices of the resolvents in the previous expression are in the set $\{N+1,\dots, N+M\}$, but this is not the case since $n\in\{2,\dots, N\}$. Hence using a fluctuation averging for $\sum_{n=2}^N\mathfrak{G}_{i_kn}^{(k)}\mathfrak{G}_{ni_2}^{(1)}$ and the one side resolvent expansion in \eqref{oneside} as in \eqref{supereqres} in the proof of Lemma \ref{selfconeq} we get \begin{align*}
-\frac{1}{\sqrt{MN}}\sum_{n=2}^N &-\left(\mathfrak{G}_{i_ki_1}^{(k)}\mathfrak{G}_{ni_1}^{(k)}+\mathfrak{G}_{i_k n}^{(k)}\mathfrak{G}_{i_1i_1}^{(k)}\right)\mathfrak{G}_{n i_2}^{(1)}\dots \mathfrak{G}_{i_{k-1}i_k}^{(k-1)}\numberthis \label{eq7}\\
&=-m_1(\zeta_1)m_1(\zeta_k)m_2(\zeta_k)\sum_{m=N+1}^{N+M} \mathfrak{G}_{mi_2}^{(1)}\dots\mathfrak{G}_{i_k m}^{(k)}+\mathcal{O}_\prec \left(\frac{1}{N^{\frac{k}{2}+\frac{1}{2}}\sqrt{\eta\eta_k}}\right).
\end{align*} 

Furthermore, following the proof of Lemma 4.1 in \cite{erdHos2017fluctuations} for the estimate of the error we obtain that \begin{equation}
\label{eq8}
R\prec \sum_{a=1}^k\frac{N^\epsilon}{\sqrt{N\eta\eta_a}}.
\end{equation} Hence, using $z_l=\zeta_l^2$ for $l=1,\dots,k$, combining \eqref{eq3} and \eqref{eq5}-\eqref{eq8} we conclude \begin{equation}
\label{aq}
\begin{split} &N^{-\frac{k}{2}}\sum_{i_1,\dots, i_k=N+1}^{N+M\, *} \mathbb{E}_{i_1}\left(\mathfrak{G}_{i_1i_2}^{(1)}\dots \mathfrak{G}_{i_{k-1}i_k}^{(k-1)}\mathfrak{G}_{i_ki_1}^{(k)}\right)\\
&\qquad =\frac{m_2(\zeta_1)}{m_1(\zeta_1)m_1(\zeta_k)m_2(\zeta_1)m_2(\zeta_k)-1}\cdot\mathcal{O}_\prec\left(\sum_{a=1}^k\frac{N^\epsilon}{\sqrt{N\eta\eta_a}}\right)\\
&\qquad=\mathcal{O}_\prec\left(\sum_{a=1}^k\frac{N^\epsilon}{(\eta_1+\eta_k)\sqrt{N\eta\eta_a}}\right),
\end{split}\end{equation} where in the last equality we used \eqref{ser} and, since \eqref{goodbound} holds true also substituting $\phi$ with $\phi^{-1}$ (see proof of Lemma \ref{sbound} in Appendix A), that $|m_2|\le \phi^{-\frac{1}{4}}\le 1$ to estimate the error. With these computations we conclude the estimate of the first term in the right-hand side of \eqref{eq3}. Notice that the estimate of the error in \eqref{aq} depends only on the Stieltjes transforms $m_1$ and $m_2$, hence, using a similar argument as in the proof of Lemma \ref{selfconeq}, we conclude that all the terms in the right-hand side of \eqref{eq3} give the same contribution. This concludes the proof of \eqref{eq1}.

The proof of \eqref{eq2}, using the equality in \eqref{aaa}, is exactly the same of \eqref{eq1} using that for the case $a=k-1$ we have the following estimate \[-\left(\mathfrak{G}_{i_{k-1}i_1}^{(k-1)}\mathfrak{G}_{ni_k}^{(k-1)}+\mathfrak{G}_{i_k n}^{(k-1)}\mathfrak{G}_{i_1i_k}^{(k-1)}\right)\mathfrak{G}_{n i_2}^{(1)}\dots \mathfrak{G}_{i_{k-2}i_{k-1}}^{(k-2)}\prec\frac{1}{N^\frac{k}{2}\sqrt{\eta\eta_{k-1}}}.\] Hence we have that \[ N^{-\frac{k+1}{2}}\sum_{i_1,\dots, i_k=N+1}^{N+M} \mathbb{E}_{i_1}\left(\mathfrak{G}_{i_1i_2}^{(1)}\dots \mathfrak{G}_{i_{k-1}i_k}^{(k-1)}\right)=\mathcal{O}_\prec\left(\sum_{a=1}^k\frac{N^\epsilon}{\sqrt{N\eta\eta_a}}\right).\] The previous expression only depends on $m_2$ and so using the same argument as before we conclude the proof of \eqref{eq2}.

The proof for $X\in\mathbb{C}^{M\times (N-1)}$ is omitted since is similar to the real case after replacing the cumulant expansion by its complex variant (Lemma 7.1 in \cite{he2017mesoscopic}).
\endproof
\end{lemma}

Notice that the estimates of Lemma \ref{cyc} hold also without the expectation:

\begin{corollary}
Under the hypotheses of Lemma \ref{cyc}, we have that for closed cycles of length $k> 2$ \begin{equation}
\label{eqAAA}
N^{-\frac{k}{2}}\sum_{j_1,\dots, j_k=1}^{M\,\, *} G_{j_1j_2}^{(1)}\dots G_{j_{k-1}j_k}^{(k-1)}G_{j_kj_1}^{(k)}\prec \frac{|z_1\dots z_k|^{-\frac{1}{2}}}{(\max_a\eta_a)\sqrt{N\eta_1\dots\eta_k}}\sum_{a=1}^k\frac{1}{\sqrt{\eta_a}},
\end{equation} and for open cycles of length $k> 1$ \begin{equation}
\label{eq9}
N^{-\frac{k+1}{2}}\sum_{j_1,\dots, j_k=1}^{M\,\, *} G_{j_1j_2}^{(1)}\dots G_{j_{k-1}j_k}^{(k-1)}\prec\frac{|z_1\dots z_k|^{-\frac{1}{2}}}{\sqrt{N\eta_1\dots\eta_{k-1}}}\sum_{a=1}^k\frac{1}{\sqrt{\eta_a}}
\end{equation}
\proof
\normalfont
First, we recall that  $\mathfrak{G}(z)$, $z\in \mathbb{C}\setminus \mathbb{R}$, is the resolvent of the linearized matrix $\mathcal{H}$. In order to prove the bounds \eqref{eqAAA}--\eqref{eq9}, we rely on \cite[Proposition~6.1]{erdHos2013averaging} with exactly the same modification as in the proof of \eqref{me}, i.e. the case when different  resolvent factors 
$\mathfrak{G}$ may have different spectral parameters. In particular, for any fixed and distinct $i_2,\dots,i_k$, 
the quantity
\begin{equation}
\label{newave}
\frac{1}{N}\sum_{i_1=N+1}^{N+M\,\, *}(1-\mathbb{E}_{i_1}) \mathfrak{G}_{i_1i_2}^{(1)}\dots \mathfrak{G}_{i_{k-1}i_k}^{(k-1)}\mathfrak{G}_{i_ki_1}^{(k)},
\end{equation}
 is smaller than the bound given by the local law of an additional factor $1/\sqrt{N\eta_1}+\dots+1/\sqrt{N\eta_k}$. Hence, the bounds in \eqref{eqAAA} and \eqref{eq9} follow by Lemma \ref{cyc}, using the relation \eqref{aaa} and that $G_{ij}=\mathcal{G}_{N+i, N+j}$ for $i,j=1,\dots,M$.
\endproof
\end{corollary}

The following lemma shows that the leading order terms of $\mathbb{E}_{1}A(z_1)\dots A(z_k)$  are the cycles of length two (see the  proof of Lemma 4.3 in \cite{erdHos2017fluctuations}).

\begin{lemma}
\label{pairs}
For each $k\ge 2$ and $z_1, \dots, z_k\in\mathbb{C}$ with $|\Im z_l|=\eta_l>0$ we have that \begin{equation}
\begin{split}
\mathbb{E}_{1}A(z_1)\dots A(z_k)&=\sum_{\pi\in P_2([k])}\prod_{\{a,b\}\in\pi}\mathbb{E}_{1}(A(z_a)A(z_b))\\
&\quad +\mathcal{O}_\prec\left(\frac{{|z_1\dots z_k|^{-\frac{1}{2}}}}{\sqrt{N\eta_1\dots\eta_k}}\sum_{a\ne b}\frac{1}{(\eta_a+\eta_b)\sqrt{\eta_a}}\right),
\end{split}\end{equation} where $[k]:=\{1,\dots,k\}$ and $P_2(L)$ is the set of pairings of the set $L$.
\end{lemma}

By Lemma \ref{pairs} we conclude that \begin{equation}
\label{final}
\begin{split}
\mathbb{E}\left[-\Im\int_\mathbb{R}\int_{\eta_0}^{10} g(z)A(z)\,d\eta dx\right]^k&=\sum_{\pi\in P_2([k])}(2V_{f,1}+(\sigma_4-1)V_{f,2})^\frac{k}{2}+\mathcal{O}_\prec\left((N^{-\frac{7}{6}}\right)\\
&=(k-1)!!(2V_{f,1}+(\sigma_4-1)V_{f,2})^\frac{k}{2}+\mathcal{O}_\prec\left(N^{-\frac{7}{6}}\right),
\end{split}\end{equation} if $k$ is even and \begin{equation}
\mathbb{E}\left[-\Im\int_\mathbb{R}\int_{\eta_0}^{10} g(z)A(z)\,d\eta dx\right]^k=\mathcal{O}_\prec\left(N^{-\frac{7}{6}}\right)
\end{equation} if $k$ is odd. If $X\in \mathbb{C}^{M\times (N-1)}$, following the same argument, we find
\[
\mathbb{E}\left[-\Im\int_\mathbb{R}\int_{\eta_0}^{10} g(z)A(z)\,d\eta dx\right]^k=(k-1)!!(V_{f,1}+|\sigma_2|^2V_{\sigma_2}+(\sigma_4-1)V_{f,2})^\frac{k}{2}+\mathcal{O}_\prec\left(N^{-\frac{7}{6}}\right).
\] 

In this way we conclude the computations of the moments for each $k\ge1$ and so with this result we have shown that the random variable $\sqrt{N}(f_N-\Omega_f)$ converges in distribution to a Gaussian random variable $\Delta_f$ with mean zero and variance $V_f$ and that any fixed moment of $\sqrt{N}(f_N-\Omega_f)$ converges to the corresponding Gaussian moment with overwhelming probability at least at a rate $\mathcal{O}\left(N^{-\frac{1}{6}+\epsilon}\right)$.

\appendix

\section{Proof of Lemma \ref{sbound}.}

We present the proof of Lemma \ref{sbound} only for $\phi\ge 1$. The proof for $\phi\le 1-d_*$ is completely analogous and so omitted.

We recall that $w_\phi(z)$ is the Stieltjes transform of the Wigner semicircle law centered in $\phi^\frac{1}{2}+\phi^{-\frac{1}{2}}$ defined as in \eqref{wigner1}. By the proof of Lemma 3.7 in \cite{erdHos2017fluctuations} and Lemma~3.6 in \cite{bloemendal2016principal}, for each $z=x+i\eta$ such that $|z-\sqrt{\phi}|\le 10$, we have that \begin{equation}
\label{gw}
c\le |w_\phi(z)|\le 1,\,\,\,|1-w_\phi(z)^2|\asymp \sqrt{\kappa_x+\eta}, \,\,\, \Im w_\phi(z)\asymp \begin{cases}
\sqrt{\kappa_x+\eta}\,\,\, \text{if} \,x\in [\gamma_-,\gamma_+]\\
\frac{\eta}{\sqrt{\kappa_x+\eta}}\,\,\,\,\,\,\,\,\, \text{if}\, x\notin [\gamma_-,\gamma_+],\end{cases}\end{equation} where $\kappa_x=\min\{|\gamma_+-x|,|\gamma_--x|\}$, $w_\phi(z):=\sqrt{\phi}(1+zm_{\phi^{-1}}(z))$ and $c>0$ is a constant independent of $\phi$.

\textbf{Proof of Lemma \ref{sbound}.} Let $\tilde{z}:=z\phi^{-\frac{1}{2}}$, taking the imaginary part of $-\frac{1}{m_\phi}=z+\tilde{z}m_\phi-(\phi^\frac{1}{2}-\phi^{-\frac{1}{2}})$ and $-\frac{1}{\tilde{z}m_\phi}=\phi^\frac{1}{2}+m_\phi-\frac{1}{\tilde{z}}(\phi^\frac{1}{2}-\phi^{-\frac{1}{2}})$ (see \eqref{st1}), we get \begin{equation}
\label{ImA}
\frac{\Im m_\phi}{|m_\phi|^2}=\eta+\Im (\tilde{z}m_\phi),\,\,\,\,\,\,\,\, \frac{\Im(\tilde{z}m_\phi)}{|\tilde{z}m_\phi|^2}=\Im m_\phi+\frac{\phi-1}{|z|^2} \eta.
\end{equation} Combining these equalities we obtain \[|\tilde{z}|^2|m_\phi|^4=1-\frac{|m_\phi|^2+\frac{\phi-1}{|z|^2}}{\Im m_\phi+\frac{\eta(\phi-1)}{|z|^2}}\eta.\] By our hypotheses $|z-\sqrt{\phi}|\le 10$ and $\phi\ge 1$, we have that $\eta\le 10$ and that there exists a constant $d>0$ independent of $\phi$ such that $|z|\le d \sqrt{\phi}$. Furthermore, from \eqref{ImA} and $\Im(\tilde{z}m_\phi)=\Im w_{\phi^{-1}}\le 1$ we have $\Im m_\phi\le C|m_\phi|^2$, with $C>0$ some constant independent of $\phi$. We conclude that 
\[|\tilde{z}|^2|m_\phi|^4=1-\frac{|m_\phi|^2+\frac{\phi-1}{|z|^2}}{\Im m_\phi+\frac{\eta(\phi-1)}{|z|^2}}\eta\le 1-2\tilde{c}\eta,\]
for any $\phi\ge 1$. The above inequality proves the bound in \eqref{goodbound}.

Furthermore, since $w_\phi(z)=-zm_\phi(z)m_{\phi^{-1}}(z)$ by \eqref{wwin} and using that, by similar computations substituting $\phi$ with $\phi^{-1}$, we have an upper bound as in \eqref{goodbound} for $|m_{\phi^{-1}}|$ and that $|w_\phi|\ge c$ from \eqref{gw}, we also obtain the lower bound in \eqref{goodbound}. Note that by a direct computation, substituting $\phi$ with $\phi^{-1}$, we get a lower bound as in \eqref{goodbound} also for $|m_{\phi^{-1}}|$. Finally, since \[1-w_\phi^2(z)=1-w_\phi(z)w_{\phi^{-1}}(z)=zm_\phi(z)+zm_{\phi^{-1}}(z)+z^2m_\phi(z)m_{\phi^{-1}}(z),\] using \eqref{win} for $zm_{\phi^{-1}}(z)$ in the right-hand side, we get \[\left|1-z\phi^{-\frac{1}{2}}m_\phi(z)^2\right|=\frac{|1-w_\phi^2(z)|}{|zm_{\phi^{-1}}(z)|}.\] Hence, using \eqref{gw} and that $|m_{\phi^{-1}}|\ge c\phi^{-\frac{1}{4}}|z|^{-\frac{1}{2}}$, we conclude \begin{equation}
\label{square}
\left|1-z\phi^{-\frac{1}{2}}m_\phi(z)^2\right|\asymp \frac{\phi^\frac{1}{4}}{|z|^\frac{1}{2}}\sqrt{\kappa_x+\eta}.
\end{equation}
This proves \eqref{finaleq}. Then, using \eqref{goodbound}, \eqref{square} and the explicit expression \[m_\phi(z)'=\frac{m_\phi(z)^2+\frac{m_\phi(z)^3}{\sqrt{\phi}}}{1-\frac{z}{\sqrt{\phi}}m_\phi(z)^2},\] obtained differentiating \eqref{st1}, we also get the bound in \eqref{boundderivative} for $|m_\phi(z)'|$.

\bibliography{Bibliography2}
\bibliographystyle{siam}
\end{document}